\documentclass[12pt]{amsart}
\usepackage{amssymb}
\usepackage{enumitem}
\usepackage{graphicx}
\makeatletter
\@namedef{subjclassname@2020}{%
	\textup{2020} Mathematics Subject Classification}
\makeatother

\usepackage[T1]{fontenc}

\newtheorem{theorem}{Theorem}[section]
\newtheorem{corollary}[theorem]{Corollary}
\newtheorem{lemma}[theorem]{Lemma}
\newtheorem{proposition}[theorem]{Proposition}

\theoremstyle{definition}
\newtheorem{definition}[theorem]{Definition}
\newtheorem{remark}[theorem]{Remark}
\newtheorem{example}[theorem]{Example}

\numberwithin{equation}{section}

\begin{document}

	\baselineskip=17pt

\author[X. Q. Qiang]{Xiangqi Qiang}
	\address{School of Mathematical Science \\ Yangzhou University\\
		Yangzhou 225002, China}
	\email{qiang.xq@qq.com}
	\author[C. J. Hou]{Chengjun Hou}
\address{School of Mathematical Science \\ Yangzhou University\\
	Yangzhou 225002, China}
\email{cjhou@yzu.edu.cn}
\thanks{Corresponding author: Chengjun Hou}	
\title[COE for automorphism systems of equivalence relations]{Continuous orbit equivalence for automorphism systems of equivalence relations} 	
	
	\begin{abstract}
	 
	We introduce notions of continuous orbit equivalence and strong (respective, weak) continuous orbit equivalence for automorphism systems of \'{e}tale equivalence relations, and characterize them in terms of the    semi-direct product groupoids, as well as their reduced groupoid $C^*$-algebras and the associated $C^*$-automorphism systems of group actions or coactions on them. In particular,  we study topological rigidity of expansive automorphism actions on compact (connected) metrizable groups.	
	\end{abstract}
	
	
	
	\maketitle

	\normalsize

		\baselineskip=17pt

\section{Introduction }

The interplay between  orbit equivalence of topological dynamical systems and $C^*$-algebras has been studied by many authors. An early celebrated result in this direction is the work on strong orbit equivalence of minimal homeomorphisms on Cantor sets given by  Giordano, Putnam and Skau (\cite{Giordano1995}). Later, Tomiyama and Boyle-Tomiyama studied a generalization of the GPS's result to the case of topologically free homeomorphisms on compact Hausdorff spaces (\cite{Boyle1998,Tomiyama1996}). In \cite{Matsumoto2010}, Matsumoto introduced the notion of continuous orbit equivalence of one-sided topological Markov shifts and characterized them in terms of  the existence of diagonals preserving $\ast$-isomorphisms between the associated Cuntz-Krieger algebras. In \cite{Matsumoto2014}, Matui and Matsumoto gave a classification result of two-sided irreducible topological Markov shifts  in the sense of flow equivalence by means of continuous orbit equivalence of one-sided topological Markov shifts. We can refer to \cite{Carlsen2019,Carlsen2017} for some generalizations on flow equivalence and study on the relation between topological conjugacy of two-sided shifts of finite type and the associated stabilized Cuntz-Krieger algebras with the canonical Cartan subalgebras and gauge actions. More recently, in \cite{Matsumoto2019,Matsumoto20211}, Mastumoto introduced notions of asymptotic continuous orbit equivalence, asymptotic conjugacy  and asymptotic flip conjugacy in Smale spaces and characterized them in terms of their groupoids and asymptotic Ruelle algebras with their dual actions, he also  characterized topological conjugacy classes of one-sided topological Markov
shifts in terms of the associated Cuntz-Krieger algebras and its gauge actions with potentials in \cite{Matsumoto20213}.

Our interests lie in group actions.  As a topological analogue of the classification results on the probability measure preserving actions in the sense of orbit equivalence, Li introduced the notion of continuous orbit equivalence for continuous group actions and proved  that two topologically free systems are continuously orbit equivalent if and only if their associated transformation groupoids are isomorphic (\cite{Li2018}). By Renault's result in \cite{Renault2008}, these conditions are also equivalent to the existence of $C^*$-isomorphism preserving the canonical Cartan subalgebras between the corresponding crossed product algebras. In \cite{Carlsen2021}, Li's rigidity result is generalized to the case of group actions with torsion-free and abelian essential stabilisers.

The local conjugacy relations from  expansive group action systems are generalizations of asymptotic equivalence relations of Smale spaces (\cite{Putnam1996,Thomsen2010}). In \cite{Hou2021}, we characterized continuous orbit equivalence of expansive systems up to local conjugacy relations and showed that two expansive actions are asymptotically continuous orbit equivalent if and only if the associated semi-direct product groupoids  of local conjugacy relations are isomorphic.

In this paper we consider continuous orbit equivalence between automorphism systems of \'{e}tale equivalence relations. Given an \'{e}tale equivalence relation $\mathcal{R}$ on a compact metrizable space $X$, let $G\curvearrowright_{\alpha} (X,\mathcal{R})$ be a dynamical system arising from an automorphism action of a countable group $G$ on $\mathcal{R}$ in the sense that each $\alpha_g$ is an automorphism of $\mathcal{R}$ as \'{e}tale groupoids. Denote by $\mathcal{R}\rtimes_{\alpha} G$ the associated semi-direct product groupoid. We say that two systems $G\curvearrowright_{\alpha} (X, \mathcal{R})$ and $H\curvearrowright_{\beta} (Y,\mathcal{S})$ \emph{conjugate} if there exist an isomorphism $\widetilde{\varphi}:\, \mathcal{R}\rightarrow \mathcal{S}$ as \'{e}tale groupoids and a group isomorphism $\theta:\, G\rightarrow H$ such that $\widetilde{\varphi}(g\gamma)=\theta(g)\widetilde{\varphi}(\gamma)$ for $\gamma\in \mathcal{R}$ and $g\in G$.  We call the set $[x]_{G,\mathcal{R}}=\{y\in X:\, (gx,y)\in \mathcal{R} \hbox{ for some $g\in G$}\}$ \emph{the bi-orbit} of $x$. Motivated by the notion of usual orbit equivalence of dynamical systems, we say that $G\curvearrowright (X,\mathcal{R})$ and $H\curvearrowright (Y,\mathcal{S})$ are \emph{orbit equivalent} if there exists a homeomorphism $\varphi:\, X\rightarrow Y$ such that $\varphi([x]_{G, \mathcal{R}})=[\varphi(x)]_{H,\mathcal{S}}$ for $x\in X$. We call they are \emph{continuously orbit equivalent} if there exist a homeomorphism $\varphi:\, X\rightarrow Y$, continuous maps $a:\,\, \mathcal{R}\times G\rightarrow H$ and  $b:\,\, \mathcal{S}\times H\rightarrow G$ such that
both the maps
$((x,y),g)\in \mathcal{R}\times G\rightarrow (\varphi(x), a((x,y),g)\varphi(g^{-1}y))\in \mathcal{S} $
and
$((x,y),g)\in \mathcal{S}\times H\rightarrow (\varphi^{-1}(x), b((x,y),g)\varphi^{-1}(g^{-1}y))\in \mathcal{R}$
are well-defined and continuous. The followings are main results in this paper.

\begin{theorem} \quad Assume that $G\curvearrowright_{\alpha} (X,\mathcal{R})$ and $H\curvearrowright_{\beta} (Y,\mathcal{S})$ are essentially free. Then the following statements are equivalent.
	\begin{enumerate}
		\item[(i)] $G\curvearrowright_{\alpha} (X,\mathcal{R})$ and $H\curvearrowright_{\beta} (Y,\mathcal{S})$ are continuously orbit equivalent;
		\item[(ii)] $\mathcal{R}\rtimes_{\alpha} G$ and $\mathcal{S}\rtimes_{\beta} H$ are isomorphic as \'{e}tale groupoids;
		\item[(iii)] there exists a $C^*$-isomorphism $\Phi$ from $C_r^*(\mathcal{R}\rtimes_{\alpha} G)$ onto $C_r^*(\mathcal{S}\rtimes_{\beta} H)$ such that $\Phi(C(X))=C(Y)$.
	\end{enumerate}
\end{theorem}

Here the notion of essential freeness for $G\curvearrowright_{\alpha} (X, \mathcal{R})$ is a generalization and analogue of topological freeness of dynamical systems. When $\mathcal{R}=\{(x,x):\, x\in X\}$ is a trivial \'{e}tale equivalence relation, or $\mathcal{R}$ is the local conjugacy relation or asymptotical equivalence relation arising from an expansive system $G\curvearrowright_{\alpha} X$ or an irreducible Smale space $(X,\varphi)$, this result is reduced to Theorem 1.2 in \cite{Li2018}, Theorem 3.4 in \cite{Matsumoto2019} and Theorem 4.2 in \cite{Hou2021}. The properties of strong or weak continuous orbit equivalence for automorphism systems are corresponding to two special orbit equivalence with some uniform conditions, and are also analogues of asymptotic flip conjugacy in \cite{Matsumoto2019} and (strong) asymptotic conjugacy in \cite{Hou2021}. Let $\rho_{\alpha}$ be the canonical cocycle from $\mathcal{R}\rtimes_{\alpha} G$ onto $G$. It follows from \cite[Lemma 6.1]{Carlsen2021} that $\rho_{\alpha}$ gives us a coaction system $(C_r^*(\mathcal{R}\rtimes_{\alpha} G), G, \delta_{\alpha})$.

\begin{theorem}\quad Assume that $G\curvearrowright_{\alpha} (X,\mathcal{R})$ and $H\curvearrowright_{\beta} (Y,\mathcal{S})$ are essentially free. Then
	\begin{enumerate}
		\item[(i)] $G\curvearrowright_{\alpha} (X,\mathcal{R})$ and $H\curvearrowright_{\beta} (Y,\mathcal{S})$ are weakly continuously orbit equivalent if and only if there is an isomorphism $\Lambda:\,\,\mathcal{R}\rtimes_{\alpha} G\rightarrow \mathcal{S}\rtimes_{\beta} H$ such that $\Lambda(\mathcal{R})=\mathcal{S}$.  Moreover, if these conditions hold, then there is a $C^*$-isomorphism $\Phi: C_r^*(\mathcal{R}\rtimes_{\alpha} G)\rightarrow C_r^*(\mathcal{S}\rtimes_{\beta} H)$ such that $\Phi(C(X))=C(Y)$ and $\Phi(C_r^*(\mathcal{R}))=C_r^*(\mathcal{S})$.
		\item[(ii)] $G\curvearrowright_{\alpha} (X,\mathcal{R})$ and $H\curvearrowright_{\beta} (Y,\mathcal{S})$ are strongly continuously orbit equivalent if and only if   there exist an \'{e}tale groupoid isomorphism $\Lambda:\, \mathcal{R}\rtimes_{\alpha} G\rightarrow \mathcal{S}\rtimes_{\beta} H$ and a group isomorphism $\theta: G\rightarrow H$ such that $\theta\rho_{\alpha}=\rho_{\beta}\Lambda$ if and only if two coaction systems $(C_r^*(\mathcal{R}\rtimes_{\alpha} G), G, \delta_{\alpha})$ and $(C_r^*(\mathcal{S}\rtimes_{\alpha} H), H, \delta_{\beta})$ are conjugate by a conjuagcy $\phi$ with $\phi(C(X))=C(Y)$.
	\end{enumerate}
	
	Furthermore,  when $\mathcal{R}$ and $\mathcal{S}$ are minimal or $X$ and $Y$ are connected, these two notions are consistent.
\end{theorem}

The assumption of essential freeness in the above theorems is necessary. Automorphism systems on local conjugacy relations from expansive actions are typical examples. The automorphism systems of local conjugacy relations from a full shift $G\curvearrowright A^G$ over a finite set $A$ and an irreducible Smale space $(X,\psi)$ are essentially free (\cite{Hou2021,Matsumoto2019}). The following result generalizes Matsumoto's result.

\begin{theorem}
	\quad
	Let $\mathcal{R}_{\alpha}$ be the local cnjugacy relation from an expansive and transitive action $G\curvearrowright_{\alpha} X$. Assume that $X$ is infinite and has no isolated points and $G$ is an abelian group such that every subgroup generated by $g$ ($g\neq e$) has finite index in $G$. Then $G \curvearrowright_{\alpha} (X,\mathcal{R}_{\alpha})$ is essentially free.
	
	Moreover, if $\mathbb{Z}\curvearrowright_{\alpha} X$ is generated by an expansive  homeomorphism $\varphi$ on $X$, then the  transitivity condition on $\varphi$ is not necessary.
\end{theorem}

In \cite{Bhattacharya2000},   Bhattacharya proved that topological conjugacy and algebraic conjugacy between two automorphism actions on compact abelian connected metrizable spaces are agreement.  We have a rigidity result for automorphism actions on nonabelian groups.

\begin{proposition} \quad Let $G\curvearrowright_{\alpha} (X,\mathcal{R})$ and $H\curvearrowright_{\beta} (Y, \mathcal{S})$ be two systems on local conjugacy relations from topologically free, expansive automorphism actions on compact and connected metrizable groups $X$ and $Y$, respectively. Assume that the homoclinic group $\Delta_{\alpha}$ associated to $G\curvearrowright_{\alpha} X$ is dense in $X$. Then the following statements are equivalent:
	\begin{enumerate}
		\item[(i)] $G\curvearrowright_{\alpha} (X,\mathcal{R})$ and $H\curvearrowright_{\beta} (Y, \mathcal{S})$ are conjugate;
		\item[(ii)] $G\curvearrowright_{\alpha} (X,\mathcal{R})$ and $H\curvearrowright_{\beta} (Y, \mathcal{S})$ are weakly continuously orbit equivalent;
		\item[(iii)] $G\curvearrowright_{\alpha} X$ and $H\curvearrowright_{\beta} Y$ are conjugate;
		\item[(iv)] $G\curvearrowright_{\alpha} X$ and $H\curvearrowright_{\beta} Y$ are algebraically conjugate.
	\end{enumerate}
	
	In particular, two hyperbolic toral automorphisms  on $\mathbb{R}^n/\mathbb{Z}^n$ are flip conjugate if and only if the $\mathbb{Z}$-actions they generates are continuously orbit equivalent up to the associated local conjugacy relations.
\end{proposition}

This paper is organized as follows. Section 3 characterizes conjugacy of automorphism systems of \'{e}tale equivalence relations and the reduced $C^*$-algebra of the associated semi-direct product groupoid of equivalence relations.  In Section 4, we introduce  notions of   continuous orbit equivalence, strong-  and weak- continuous orbit equivalence for automorphism systems, and characterize them in terms of the semi-direct product groupoids and the corresponding $C^*$-algebras. In Section 5, we discuss essential freeness of automorphism systems on local conjugacy equivalence relations arising from expansive actions, and in Section 6, we study topological rigidity of expansive automorphism actions on compact (connected) metrizable groups. As an example, we characterize the structure of the local conjugacy relation from a hyperbolic toral automorphism on $n$-torus.

\section{Preliminaries}

Unless otherwise specified, all our groups are discrete and countable, their identity elements are denoted by the same symbol $e$, and all topological groupoids are second countable, locally compact and Hausdorff. We refer to \cite{Renault1980,Sims2020} for more details on topological groupoids and their $C^*$-algebras, and refer to \cite{Pedersen1979,Williams2007} for $C^{*}$-dynamical systems.

For a topological groupoid $\mathcal{G}$, let $\mathcal{G}^{(0)}$ and $\mathcal{G}^{(2)}$ be  the unit space and the set of composable pairs, respectively. The range and domain maps $r,d$ from $\mathcal{G}$ onto $\mathcal{G}^{(0)}$  are defined by $r(g)=gg^{-1}$ and $d(g)=g^{-1}g$, respectively. If $r$ and $d$ are local homeomorphisms then $\mathcal{G}$ is called to be \emph{\'{e}tale}.
For $u,v\in \mathcal{G}^{(0)}$, we write $\mathcal{G}^u=r^{-1}(u)$, $\mathcal{G}_u=d^{-1}(u)$ and $\mathcal{G}_u^v=\mathcal{G}^v\cap\mathcal{G}_u$. When $\mathcal{G}$ is \'{e}tale, these sets are discrete and countable, and $\mathcal{G}^{(0)}$ is open and closed in $\mathcal{G}$. Recall that $\mathcal{G}$ is \emph{topologically principle} if $\left\lbrace u\in \mathcal{G}^{(0)}:\,\mathcal{G}_{u}^{u}=\left\lbrace u\right\rbrace \right\rbrace $ is dense in $\mathcal{G}^{(0)}$.

Each equivalence relation $\mathcal{R}\subseteq X\times X$ on a topological space $X$ is a groupoid with  multiplication $(x,y)(w,z)=(x,z)$ if $y=w$ and inverse
$(x,y)^{-1}=(y,x)$. If we identify $(x,x)$ with $x$, then the unit space $\mathcal{R}^{(0)}$ equals $X$ and the range (resp. domain) map is defined by $r(x,y)=x$ (resp. $d(x,y)=y$).
If there exists a topology on $\mathcal{R}$ (not necessarily the relative product topology from $X\times X$) for which $\mathcal{R}$ is an \'{e}tale groupoid, then $\mathcal{R}$ is called an \emph{\'{e}tale equivalence relation} on $X$. In this case, if every $\mathcal{R}$-equivalence class is dense in $X$ then $\mathcal{R}$ is \emph{minimal}.

By a \emph{dynamical system}, denoted by $G\curvearrowright_{\alpha} X$ (or simply by $G\curvearrowright X$), we mean an action $\alpha$ of a group $G$ on a second countable, locally compact and Hausdorff space $X$ by homeomorphisms. The action $\alpha$ is usually expressed as $(g,x)\in G\times X\rightarrow gx\in X$. The associated transformation groupoid $X\rtimes G$ is given by the set $X\times G$ with the product topology,  multiplication $(x,g)(y,h)=(x,gh)$ if $y=g^{-1}x$, and inverse $(x,g)^{-1}=(g^{-1}x,g^{-1})$. Clearly, $X\rtimes G$ is \'{e}tale, and if $(x,e)$ is identified with $x$ then  its unit space equals $X$, range map $r(x,g)=x$ and domain map  $d(x,g)=g^{-1}x$.
A system $G\curvearrowright X$ is  said to be \emph{topologically free} if for every $e\neq g\in G$, $\{x\in X:\, gx\neq x\}$ is dense in $X$. From \cite[Corollary 2.3]{Li2018}, $G\curvearrowright X$ is topologically free if and only if $X\rtimes G$ is topologically principal. Two systems $G\curvearrowright X$ and $H\curvearrowright Y$ are \emph{ conjugate} if there are a homeomorphism $\varphi:\, X\rightarrow Y$ and a group isomorphism $\theta:\, G\rightarrow H$ such that $\varphi(gx)=\theta(g)\varphi(x)$ for $x\in X$ and $g\in G$.

A map $\Phi:\, \mathcal{G}\rightarrow \mathcal{H}$ between \'{e}tale groupoids $\mathcal{G}$ and $\mathcal{H}$ is a homomorphism if it is continuous and, for all $(\gamma,\gamma')\in \mathcal{G}^{(2)}$,  we have $(\Phi(\gamma),\Phi(\gamma'))\in \mathcal{H}^{(2)}$ and $\Phi(\gamma\gamma')=\Phi(\gamma)\Phi(\gamma')$. Moreover, if $\Phi$ is a homeomorphism such that $\Phi$ and $\Phi^{-1}$ are homomorphisms, then it is called an \emph{ isomorphism}. In this case, the restriction, $\Phi|_{\mathcal{G}^{(0)}}$, of $\Phi$ to the unit space $\mathcal{G}^{(0)}$ is a homeomorphism from $\mathcal{G}^{(0)}$ onto $\mathcal{H}^{(0)}$. A homomorphism from  $\mathcal{G}$ into  a group $\Gamma$ is also called a \emph{cocycle} on $\mathcal{G}$.
Two  \'{e}tale equivalence relations $\mathcal{R}\subseteq X\times X$ and $\mathcal{S}\subseteq Y\times Y$ are isomorphic if and only if there exists a homeomorphism $\varphi: X \rightarrow  Y$ such
that $\varphi\times \varphi: (x,y)\in\mathcal{R}\rightarrow (\varphi(x),\varphi(y))\in \mathcal{S}$ is an isomorphism.

Given an \'{e}tale groupoid $\mathcal{G}$, the linear space, $C_c(\mathcal{G})$, of continuous complex functions  with compact support
on $\mathcal{G}$ is a $\ast$-algebra under the operations: $f^{\ast}(\gamma)=\overline{f(\gamma^{-1})}$ and $f\ast
g(\gamma)=\sum_{\gamma'\in \mathcal{G}_{d(\gamma)}}f(\gamma\gamma'^{-1})g(\gamma')$ for $f,g\in C_c(\mathcal{G})$ and $\gamma\in \mathcal{G}$.
For each $u\in \mathcal{G}^{(0)}$, there is a $\ast$-representation $Ind_{u}$ of $C_c(\mathcal{G})$ on the Hilbert space $l^2(\mathcal{G}_u)$ of square summable functions on $\mathcal{G}_u$  by
$Ind_{u}(f)(\xi)(\gamma)=\sum_{\gamma'\in \mathcal{G}_u}f(\gamma\gamma'^{-1})\xi(\gamma')$  for $f\in C_c(\mathcal{G})$, $\xi\in l^2(\mathcal{G}_u)$ and $\gamma\in \mathcal{G}_u$. The reduced $C^*$-algebra $C_r^*(\mathcal{G})$ of $\mathcal{G}$ is the completion of $C_c(\mathcal{G})$ with respect to the norm $\|f\|_{red}=\sup_{u\in \mathcal{G}^{(0)}}\|Ind_{u}(f)\|$ for $f\in C_c(\mathcal{G})$.  Since $\mathcal{G}^{(0)}$ is clopen in $\mathcal{G}$, $C_c(\mathcal{G}^{(0)})$ is contained in $C_c(\mathcal{G})$ in the canonical way, and this extends to an injection $C_0(\mathcal{G}^{(0)})\hookrightarrow C_r^*(\mathcal{G})$. For an open subgroupoid $\mathcal{H}$ of $\mathcal{G}$, $C_c(\mathcal{H})$ can be embedded into $C_c(\mathcal{G})$ as a $\ast$-subalgebra, so $C_r^*(\mathcal{H})$ is embedded into $C_r^*(\mathcal{G})$ as a $C^*$-subalgebra in the canonical way. The $C^*$-algebra $C_r^*(X\rtimes G)$ of the transformation groupoid is isomorphic to the reduced crossed product $C_0(X)\times_{\alpha, r} G$ (\cite{Sims2020}).

Given two groups $N$, $H$ and a homomorphism $\varphi$ from $H$ into the automorphism group $Aut(N)$ of $N$, the semi-direct product, denoted by $N\rtimes_{\varphi}H$, of $N$ by $H$ is defined as the set $N\times H$ with group law given by the formulas $(n,h)(n_1,h_1)=(n\varphi_{h}(n_1),hh_1)$ and
$(n,h)^{-1}=(\varphi_{h^{-1}}(n^{-1}),h^{-1})$.

\section{Automorphism systems of \'{e}tale equivalence relations and the associated semi-direct product groupoids}

Given an \'{e}tale equivalence relation $\mathcal{R}$ on a compact metrizable space $X$, we call a dynamical system $G\curvearrowright_{\alpha} \mathcal{R}$ an \emph{automorphism system} if each $\alpha_g$ is an automorphism of $\mathcal{R}$ as \'{e}tale groupoids. Clearly, this system induces an action, also denoted by $\alpha$, of $G$ on $X$ by homeomorphisms such that $g(x,y)=(gx,gy)$ for $g\in G$ and $(x,y)\in \mathcal{R}$. We use the notation $G\curvearrowright_{\alpha} (X,\mathcal{R})$ (or $G\curvearrowright (X,\mathcal{R})$ for short) to denote such automorphism system.

The semi-direct product groupoid, $\mathcal{R}\times_{\alpha} G$, attached to $G\curvearrowright_{\alpha} (X,\mathcal{R})$,  is the set $\mathcal{R}\times G$ with inverse $((x,y),g)^{-1}=((g^{-1}y,g^{-1}x), g^{-1})$, and multiplication $((x,y),g)((u,v),h)=((x,gv),gh)$ if $u=g^{-1}y$. The unit space identifies with $X$ by identifying $((x,x),e)$ with $x$. Then $r((x,y),g)=x$ and $d((x,y),g)=g^{-1}y$. Endowed with the relative product topology from $\mathcal{R}\times G$,  the groupoid $\mathcal{R}\times_{\alpha} G$ is \'{e}tale (\cite{Renault1980}). The following is another characterization of the semi-direct product groupoid.

\begin{definition}\quad Let
	$$\mathcal{R}\rtimes_{\alpha} G=\left\lbrace (x,g,y) | \ g\in G,x,y\in X,(x,gy)\in \mathcal{R}\right\rbrace. $$
	Then, under the following multiplication and inverse,
	$$(x,g,y)(y,h,v)=(x,gh,v), \,\, \hbox{ and } (x,g,y)^{-1}=(y,g^{-1},x),$$
	$\mathcal{R}\rtimes_{\alpha} G$ is a groupoid. Define a map $\gamma_{0}:\, \mathcal{R}\rtimes_{\alpha} G \rightarrow \mathcal{R}\times_{\alpha} G$, by
	$\gamma_{0}(x,g,y)=((x,gy),g)$, which is a bijection with inverse $\gamma_{0}^{-1}((x,y),g)=(x,g,g^{-1}y)$. We transfer the product topology from $\mathcal{R}\times_{\alpha} G$ over $\mathcal{R}\rtimes_{\alpha} G$. Then
	$\mathcal{R}\rtimes_{\alpha} G$ is an \'{e}tale groupoid and $\gamma_0$ is an \'{e}tale groupoid isomorphism.
\end{definition}

\begin{remark}\quad
	If we identify the unit space $(\mathcal{R}\rtimes_{\alpha} G)^{(0)}$ with $X$ as topological spaces by identifying $(x,e,x)$ with $x$, then $r(x,g,y)=x$ and $d(x,g,y)=y$. The equivalence relation $\mathcal{R}$ and the transformation groupoid $X\rtimes G$ can be embedded into $\mathcal{R}\rtimes_{\alpha} G$ as \'{e}tale subgroupoids  through the identifications
	$(x,y)\in \mathcal{R}\rightarrow (x,e,y)\in \mathcal{R}\rtimes_{\alpha} G$ and $(x,g)\in X\rtimes G\rightarrow (x,g,g^{-1}x)\in \mathcal{R}\rtimes_{\alpha} G$.
	
	One can check that the map  $\rho_{\alpha}:\, \mathcal{R}\rtimes_{\alpha} G\rightarrow G$, defined by $\rho_{\alpha}(x,g,y)=g$, is a cocycle.
\end{remark}

We call two automorphism systems $G\curvearrowright_{\alpha} (X, \mathcal{R})$ and $H\curvearrowright_{\beta} (Y,\mathcal{S})$ on compact metrizable spaces \emph{conjugate} if there are an isomorphism $\widetilde{\varphi}:\, \mathcal{R}\rightarrow \mathcal{S}$ and a group isomorphism $\theta:\, G\rightarrow H$ such that $\widetilde{\varphi}(g\gamma)=\theta(g)\widetilde{\varphi}(\gamma)$ for $\gamma\in \mathcal{R}$ and $g\in G$. Clearly, this is equivalent to that there are a homeomorphism $\varphi:\, X\rightarrow Y$ and a group isomorphism $\theta:\, G\rightarrow H$ such that $\varphi\times\varphi:\, (x,y)\in\mathcal{R}\rightarrow (\varphi(x),\varphi(y))\in\mathcal{S}$ is an isomorphism and $\varphi(gx)=\theta(g)\varphi(x)$ for $x\in X$ and $g\in G$. In particular, two systems $G\curvearrowright_{\alpha} X$ and $H\curvearrowright_{\beta} Y$ are conjugate.

\begin{proposition}\quad If $G\curvearrowright_{\alpha} (X, \mathcal{R})$ and $H\curvearrowright_{\beta} (Y,\mathcal{S})$ are conjugate, then there is an isomorphism, $\Lambda:\,\,\mathcal{R}\rtimes_{\alpha} G\rightarrow \mathcal{S}\rtimes_{\beta} H$, such that $\Lambda(\mathcal{R})=\mathcal{S}$ and $\Lambda(X\rtimes G)=Y\rtimes H$.
	
	Assume that one of the following statements holds:
	\begin{enumerate}
		\item[(i)] $X$ and $Y$ are connected.
		\item[(ii)] $\mathcal{R}$ and $\mathcal{S}$ are minimal.
	\end{enumerate}
	Then  the above converse holds, i.e., $G\curvearrowright_{\alpha} (X, \mathcal{R})$ and $H\curvearrowright_{\beta} (Y,\mathcal{S})$ are conjugate if and only if there is an isomorphism, $\Lambda:\,\,\mathcal{R}\rtimes_{\alpha} G\rightarrow \mathcal{S}\rtimes_{\beta} H$, such that $\Lambda(\mathcal{R})=\mathcal{S}$ and $\Lambda(X\rtimes G)=Y\rtimes H$.
\end{proposition}		

\begin{proof}\quad Assume that $G\curvearrowright_{\alpha} (X, \mathcal{R})$ and $H\curvearrowright_{\beta} (Y,\mathcal{S})$ are conjugate by a homeomorphism $\varphi$ from $X$ onto $Y$ and a group isomorphism $\theta$ from $G$ onto $H$. Define the map $\Lambda$ from $\mathcal{R}\rtimes_{\alpha} G$ into $\mathcal{S}\rtimes_{\beta} H$ by $\Lambda(x,g,y)=(\varphi(x),\theta(g),\varphi(y))$. Then $\Lambda$ is an isomorphism with inverse $\Lambda^{-1}(u,h,v)=(\varphi^{-1}(u),\theta^{-1}(h),\varphi^{-1}(v))$ and $\Lambda(\mathcal{R})=\mathcal{S}$ and $\Lambda(X\rtimes G)=Y\rtimes H$.
	
	On the contrary, let $\Lambda$ be an isomorphism from $\mathcal{R}\rtimes_{\alpha} G$ onto $\mathcal{S}\rtimes_{\beta} H$ such that $\Lambda(\mathcal{R})=\mathcal{S}$ and $\Lambda(X\rtimes G)=Y\rtimes H$. Let $\varphi$ be the restriction of $\Lambda$ to $X$, and let $a=\rho_{\beta}\Lambda$ and $b=\rho_{\alpha}\Lambda^{-1}$. Then $\varphi$ is a homeomorphism from $X$ onto $Y$, and $a$ and $b$ are continuous cocycles on $\mathcal{R}\rtimes_{\alpha} G$ and $\mathcal{S}\rtimes_{\beta} H$, respectively. Moreover, $\Lambda(x,g,y)=(\varphi(x), a(x,g,y),\varphi(y))$, and its inverse $\Lambda^{-1}(u,h,v)=(\varphi^{-1}(u),b(u,h,v),\varphi^{-1}(v))$. The fact that $\Lambda(\mathcal{R})=\mathcal{S}$ implies that $a(x,e,y)=e$ and $\varphi\times\varphi:\, (x,y)\in \mathcal{R}\rightarrow (\varphi(x),\varphi(y))\in \mathcal{S}$ is an isomorphism. The requirement that $\Lambda(X\rtimes G)=Y\rtimes H$ gives us that $$\varphi(x)=a(x,g,g^{-1}x)\varphi(g^{-1}x). \eqno{(3.1)}$$  Also since $(x,g,g^{-1}x)(g^{-1}x,e,g^{-1}y)(g^{-1}y,g^{-1},y)=(x,e,y)$ for $(x,y)\in \mathcal{R}$ and $g\in G$, we have $a(x,g,g^{-1}x)=a(y,g,g^{-1}y)$. By symmetry, $b$ has a similar property to $a$.

	Assume that $X$ and $Y$ are connected. Since the  restriction map $a|_{X\rtimes G}:\, X\rtimes G\rightarrow H$ is continuous, we have, for every $g\in G$, the restriction map $a|_{X\times\{g\}}$ is a constant, and thus $a(x,g,g^{-1}x)=a(y,g, g^{-1}y)$ for all $x,y\in X$ and $g\in G$.  Similarly, we have $b(u,h,h^{-1}u)=b(v,h,h^{-1}v)$ for  all $u,v\in Y$ and $h\in H$.
	
	Assume that $\mathcal{R}$ and $\mathcal{S}$ are minimal. For $x,y\in X$ and $g\in G$, we choose a sequence $\{x_n\}$ in $X$ converging to $y$ and satisfying $(x_n,x)\in \mathcal{R}$ for each $n$. From the above proof, $a(x_n,g,g^{-1}x_n)=a(x,g,g^{-1}x)$ for each $n$, which implies that $a(x,g,g^{-1}x)=a(y,g,g^{-1}y)$ from the continuity of $a$. Similarly, we have $b(u,h,h^{-1}u)=b(v,h,h^{-1}v)$ for  all $u,v\in Y$ and $h\in H$.

	Consequently, under the hypothesis of (i) or (ii), there exist two maps $\theta:\,G\rightarrow H$ and $\vartheta:\,H\rightarrow G$ such that $a(x, g,g^{-1}x)=\theta(g)$ and  $b(u, h,h^{-1}u)=\vartheta(h)$ for every $x\in X, u\in Y$, $g\in G$ and $h\in H$. Since $\Lambda$ is an isomorphism with inverse $\Lambda^{-1}$, $\theta$ is a group isomorphism with inverse $\vartheta$. Moreover, (3.1) implies that $\varphi(gx)=\theta(g)\varphi(x)$ for $x\in X$ and $g\in G$. Hence $G\curvearrowright_{\alpha} (X, \mathcal{R})$ and $H\curvearrowright_{\beta} (Y,\mathcal{S})$ are conjugate.

\end{proof}

Given an automorphism system $G\curvearrowright_{\alpha} (X, \mathcal{R})$, one can check that the map $$\alpha_g(f)(x,y)=f(g^{-1}x,g^{-1}y)$$ for $f\in C_c(\mathcal{R})$, $(x,y)\in \mathcal{R}$ and $g\in G$ gives an $C^{*}$-dynamical system $(C_{r}^{*}(\mathcal{R}),G,\alpha)$. Let $C_{c}(G,C_{r}^{*}(\mathcal{R}))$  be the  set of all continuous complex functions from $G$ to $ C^{*}_{r}(\mathcal{R})$ with compact support sets, which is a $\ast$-algebra over $\mathbb{C}$ under the  following multiplicative and convolution:
$$(\xi\ast\eta)(g)=\sum_{h\in G}\xi(h)\alpha_{h}\left( \eta\left( h^{-1}g\right) \right) $$
$$ \xi^{*}(g)=\alpha_{g}\left( \xi\left( g^{-1}\right) ^{*}\right) $$
for $\xi,\eta \in C_{c}(G,C_{r}^{*}(\mathcal{R}))$ and whose closure under the reduced crossed norm is denoted by $C_{r}^{*}(\mathcal{R}) \rtimes_{\alpha,r} G$, and is called the reduced crossed product $C^*$-algebra associated to $(C_{r}^{*}(\mathcal{R}),G,\alpha)$ (\cite{Pedersen1979,Williams2007}). By identifying an element $a\in C_r^*(\mathcal{R})$ with the element $\xi_a\in C_c(G, C_r^*(\mathcal{R})$ defined by $\xi_a(e)=a$ and $\xi_a(g)=0$ for $g\neq e$, $C_r^*(\mathcal{R})$ can be embedded into $C_r^*(\mathcal{R})\rtimes_{a,r}G$ as a unital $C^*$-subalgebra. When $G$ is abelian, we let $(C_{r}^{*}(\mathcal{R}) \rtimes_{\alpha,r} G, \widehat{G}, \widehat{\alpha})$ be the (dual) $C^*$-automorphism system of the dual group $\widehat{G}$, defined by 
$$ \widehat{\alpha}_{\xi}(f)(g)=<\xi,g>f(g)$$
for $ \xi \in \widehat{G}$, $f\in C_{c}(G,C_{r}^{*}(\mathcal{R}))$ and $g\in G$,   where $<\xi,g>$ is the value of
the character $\xi\in \widehat{G}$ at $g\in G$.

Recall that a conjugacy between two $C^*$-dynamical systems $(\mathcal{A}, G, \alpha)$ and $(\mathcal{B},H,\beta)$ is a $*$-isomorphism $\phi: \mathcal{A}\rightarrow \mathcal{B}$ that is $\alpha$ - $\beta$ equivariant in the sense that there exists  a group isomorphism $\theta: G\rightarrow H$ satisfying that $\phi\alpha_g=\beta_{\theta(g)}\phi$ for each $g\in G$. If such a $\phi$ exists, we call two systems conjugate.  Note that the existence of an isomorphism between two \'{e}tale equivalence relations $\mathcal{R}$ on $X$ and $\mathcal{S}$ on $Y$ is consistent with the existence of a $C^*$-isomorphism between their associated reduced groupoid $C^*$-algebras $C_r^*(\mathcal{R})$ and $C_r^*(\mathcal{S})$ preserving the canonical subalgebras $C(X)$ and $C(Y)$ (\cite{Renault1980}). Thus one can check the following proposition by definitions.

\begin{proposition}\quad
	Two automorphism systems $G\curvearrowright_{\alpha} (X, \mathcal{R})$ and $H\curvearrowright_{\beta} (Y,\mathcal{S})$ are conjugate if and only if there is a conjuacy $\phi$ between $(C_r^*(\mathcal{R}),G,\alpha)$ and $(C_r^*(\mathcal{S}),H,\beta)$ such that $\phi(C(X))=C(Y)$.
	
	In this case, there exists a $\ast$-isomorphism $\Lambda: C_r^*(\mathcal{R})\rtimes_{\alpha,r}G\rightarrow C_r^*(\mathcal{S})\rtimes_{\beta,r}H$ such that $\Lambda(C_r^*(\mathcal{R}))=C_r^*(\mathcal{S})$ and $\Lambda(C(X))=C(Y)$.
\end{proposition}

From \cite[Proposition II.5.1]{Renault1980}, when $G$ is abelian, the canonical cocycle $\rho_{\alpha}$ on $\mathcal{R}\rtimes_{\alpha} G$ induces the dual action, denoted by $\widehat{\rho_{\alpha}}$, of the dual group $\widehat{G}$ on $C_r^*(\mathcal{R}\rtimes_{\alpha} G)$, defined by
$$ (\widehat{\rho_{\alpha}}_{\xi}f)(x,g,y)=<\xi,g>f(x,g,y)$$
for $\xi\in \widehat{G}$, $f\in C_c(\mathcal{R}\rtimes_{\alpha} G)$, $(x,g,y)\in \mathcal{R}\rtimes_{\alpha} G$. Thus, this forms a $C^*$-dynamical system $(C_{r}^{*}( \mathcal{R}\rtimes_{\alpha} G), \widehat{G},\widehat{\rho_{\alpha}})$. Moreover, if $G=\mathbb{Z}$, then the fixed point algebra of $\widehat{\rho_{\alpha}}$ is isomorphic to $C_r^*(\mathcal{R})$ (\cite[Proposition 3.3.7]{Renault2009}).
The following theorem characterizes the reduced groupoid $C^*$-algebra of $\mathcal{R}\rtimes_{\alpha} G$ by crossed product construction, which is perhaps a well-known fact, as we were unable to find an explicit reference, we provide a proof.

\begin{theorem}\quad
	Let $G\curvearrowright_{\alpha} (X, \mathcal{R})$ be an automorphism system. Then $C_{r}^{*}( \mathcal{R}\rtimes_{\alpha} G)$ is isomorphic to $C_{r}^{*}(\mathcal{R}) \rtimes_{\alpha,r} G$.
	Moreover, if $G$ is abelian, then two $C^*$-dynamical systems $(C_{r}^{*}( \mathcal{R}\rtimes_{\alpha} G), \widehat{G},\widehat{\rho_{\alpha}})$ and $(C_{r}^{*}(\mathcal{R}) \rtimes_{\alpha,r} G, \widehat{G}, \widehat{\alpha})$ are conjugate.
\end{theorem}

\begin{proof}\quad To simplify symbol, let $\mathcal{G}=\mathcal{R}\rtimes_{\alpha} G$.  Define
	$$ \Phi(\xi)(x,g,y)=\xi(g)(x,gy)\,\hbox{ for $\xi\in  C_{c}(G,C_c(\mathcal{R}))$ and $(x,g,y)\in  \mathcal{G}$ }$$
	and $$ \Psi (\eta)(g)(x,y)=\eta(x,g,g^{-1}y),\, \hbox{ for $\eta\in C_{c}(\mathcal{G})$ and $g\in G,(x,y)\in \mathcal{R}$. }$$
	One can check that $\Phi: C_{c}(G,C_c(\mathcal{R})) \rightarrow  C_{c}(\mathcal{G})$ and $\Psi: C_{c}(\mathcal{G}) \rightarrow   C_{c}(G,C_c(\mathcal{R}))$ are $\ast$-isomorphisms such that $\Phi$ and $\Psi$ are inverse to each other.
	
	Given $x\in X$, let $l^2(\mathcal{R}_x)$ be the Hilbert space of all square-summable complex-valued functions on the $\mathcal{R}$-equivalent class $\mathcal{R}_x$ of $x$. We consider two Hilbert spaces $l^2(G,l^2(\mathcal{R}_{x}))=\{\varphi:\,\, G\rightarrow l^2(\mathcal{R}_{x})|\,\,\, \sum\limits_{g \in G}\|\varphi(g)\|^2<+\infty\}$ and
	$l^2(\mathcal{G}_{x})=\{ \psi:\,\, \mathcal{G}_{x} \rightarrow \mathbb{C}\ | \sum\limits_{\gamma \in \mathcal{G}_{x}}\|\psi(\gamma)\|^2<+\infty\}$. Then the map $U_{x}$, defined by $(U_x\varphi)(y, g, x)=\varphi(g)\left(g^{-1} y ,x\right)$ for $\varphi\in l^{2}\left(G, l^{2}\left(\mathcal{R}_{x}\right) \right)$ and $(y,g,x)\in \mathcal{G}_{x}$, is a unitary operator from $l^{2}\left(G, l^{2}\left(\mathcal{R}_{x}\right)\right)$  onto $l^{2}\left(\mathcal{G}_{x}\right)$.
	
	Let $\pi_{x}$ and $\lambda _{x}$ be the regular representations of $C_{c}(\mathcal{G})$ on $l^{2}(\mathcal{G}_{x})$ and $C_c(\mathcal{R})$ on $l^{2}(\mathcal{R} _{x})$ associated to $x$, respectively.
	Then we have the direct sums of representations
	$$\pi={\textstyle \bigoplus\limits_{x\in X}}\pi_x:\, C_c(\mathcal{G})\rightarrow {\textstyle \bigoplus\limits_{x\in X}}B(l^2(\mathcal{G}_x)),\,\,\, \lambda={\textstyle \bigoplus\limits_{x\in X}}\lambda_x:\, C_c(\mathcal{R})\rightarrow {\textstyle \bigoplus\limits_{x\in X}}B(l^2({\mathcal{R}}_x))$$
	Then $\pi_x$, $\lambda_x$, $\pi$ and $\lambda$ can be extended to their corresponding reduced groupoid $C^*$-algebras and we use the same symbols to denote their extensions. Moreover, $\pi$ and $\lambda$ are faithful representations on $C_r^*(\mathcal{G})$ and $C_r^*(\mathcal{R})$, respectively.
	
	The representation $\lambda$ induces a faithful representation $$\widetilde{\lambda}:\,\,\xi\in C_c(G,C_r^*(\mathcal{R}))\rightarrow {\textstyle \bigoplus\limits_{x\in X}}\widetilde{\lambda}_x(\xi)\in {\textstyle \bigoplus\limits_{x\in X}}B( l^2(G,l^2({\mathcal{R}}_x))),$$
	where, for each $x\in X$, $\tilde{\lambda}_{x}$ is the representation of $C_c(G,C_r^*(\mathcal{R}))$ on the Hilbert space $l^2(G,l^2({\mathcal{R}}_x))$, given by
	$(\tilde{\lambda}_{x}(\xi)\varphi) (g)= \underset{h \in G}{\sum}\lambda _{x}(\alpha_{g^{-1}}(\xi(h))) \varphi\left(h^{-1} g\right)$ for $\xi\in C_{c}( G, C_{r}^{*}(\mathcal{R}))$,
	$\varphi \in l^{2}\left(G, l^{2}\left(\mathcal{R}_{x}\right)\right)$ . Let $\widehat{\lambda}_x(\xi)= U_x\widetilde{\lambda}_x(\xi)U_x^*$ for $x\in X$ and $\xi\in C_{c}( G, C_{r}^{*}(\mathcal{R}))$. Then $$\widehat{\lambda}:\, \xi\in C_c(G,C_r^*(\mathcal{R}))\rightarrow {\textstyle \bigoplus\limits_{x\in X}} \widehat{\lambda}_x(\xi)\in {\textstyle \bigoplus\limits_{x\in X}}B( l^2(\mathcal{G}_x)),$$ is a faithful representation.
	We can check that $\pi_x\Phi(\xi)=\widehat{\lambda}_x(\xi)$ for each $x\in X$, thus $\pi\Phi(\xi)=\widehat{\lambda}(\xi)$ for all $\xi\in C_c(G,C_c(\mathcal{R}))$.
	
	In fact, for each $\varphi $ in $ {l}^{2}\left(G, l^{2}\left(\mathcal{R}_{x}\right)\right),  (y, g, x) $ in $ \mathcal{G}_{x}$, we have
	$$\begin{aligned}
		(\pi_{x}\Phi(\xi)U_{x}) &(\varphi)(y , g, x)\\ &=\sum_{h\in G,(u,h^{-1}gx) \in \mathcal{R}} [\Phi(\xi)(y, h, u)] [U_{x} \varphi\left(u, h^{-1} g, x\right)]  \\
		&=\sum _{h\in G,(g^{-1} h u, x)\in \mathcal{R}_{x} }[\xi(h)(y, h u)][\varphi\left(h^{-1} g\right)(g^{-1} hu,x)]\\
		&=\sum_{h\in G, (v,x)\in \mathcal{R}_{x}}\xi(h)(y, g v) \varphi \left(h^{-1} g\right)(v, x)
	\end{aligned}$$
	and
	$$\begin{aligned}
		U_{x}\left(\tilde{\lambda_{x}} (\xi)(\varphi)\right)(y, g, x) &=(\tilde{\lambda_{x}}(\xi)\varphi) (g)\left(g^{-1} y, x\right) \\ &=\sum _{h\in G} \sum _{(u, x) \in \mathcal{R}} \xi(h)(y, g u) \varphi\left(h^{-1} g\right)(u, x).
	\end{aligned}$$
	Then, for each $ \xi\in C_{c} ( G, C_c(\mathcal{R})) $, we have
	$$\|\Phi(\xi)\|_{r e d} =\sup _{x \in X}\|\pi_{x}(\Phi(\xi))\|_{B\left(l^{2}\left(\mathcal{G}_{x}\right)\right)} =\sup _{x \in X}\|\hat{\lambda}_{x}(\xi)\|_{\left.B ( l^{2}\left(\mathcal{G}_{x}\right)\right)}. $$
	Thus $\|\Phi(\xi)\|_{red}=\|\xi\|_{red}$ for $\xi\in C_c(G,C_c(\mathcal{R}))$, and $\Phi$ is an isomorphism.
	
	The conjugacy of two $C^*$-systems follows from the definitions of dual actions and the construction of $\Psi$. 	
\end{proof}

\begin{remark} For a countable discrete group $\Gamma$, let $\lambda:\, g\in \Gamma \rightarrow \lambda_g\in B(l^2(\Gamma))$ be the left regular representation of $\Gamma$, and $C_r^*(\Gamma)$ be the reduced group $C^*$-algebra of $\Gamma$. Let $\delta_{\Gamma}: C_r^*(\Gamma)\rightarrow C_r^*(\Gamma)\otimes C_r^*(\Gamma)$ (where we use the minimal tensor product) be the $C^*$-homomorphism defined by $\delta_{\Gamma}(\lambda_g)=\lambda_g\otimes\lambda_g$ for each $g\in \Gamma$. Given a unital $C^*$-algebra $\mathcal{A}$, we recall that a coaction of $\Gamma$ on $\mathcal{A}$ is a nondegenerate homomorphism $\delta:\, \mathcal{A}\rightarrow \mathcal{A}\otimes C_r^*(\Gamma)$ satisfying the coaction identity $(\delta\otimes id)\circ\delta=(id\otimes \delta_{\Gamma})\circ\delta$, where $id$ is the identity map. We call $(A,\Gamma,\delta)$ a $C^*$-coaction system. Recall that two $C^*$-coaction systems $(\mathcal{A},G,\delta)$ and $(\mathcal{B},H,\varrho)$ are called conjugate if there exists a conjugacy $\phi$ between two systems, that is,  $\phi$ is a $C^*$-isomorphism from $\mathcal{A}$ onto $\mathcal{B}$ such that there exists  an isomorphism $\theta: G\rightarrow H$ satisfying $(\phi\otimes \widetilde{\theta})\circ \delta=\varrho\circ \phi$, where $\widetilde{\theta}:\, C_r^*(G)\rightarrow C_r^*(H)$ is the $C^*$-isomorphism induced by $\theta$.
	
	For an automorphism system $G\curvearrowright_{\alpha} (X, \mathcal{R})$, it follows from \cite[Lemma 6.1]{Carlsen2021} that the canonical cocycle $\rho_{\alpha}: \mathcal{R}\rtimes_{\alpha} G\rightarrow G$ induces a coaction $\delta_{\alpha}:\, C_r^*(\mathcal{R}\rtimes_{\alpha} G)\rightarrow C_r^*(\mathcal{R}\rtimes_{\alpha} G)\otimes C_r^*(G)$, of $G$ on $C_r^*(\mathcal{R}\rtimes_{\alpha} G)$ such that $\delta_{\alpha}(f)=f\otimes \lambda_g$ when $g\in G$ and $f\in C_c(\mathcal{R}\rtimes_{\alpha} G)$ satisfy that $supp(f)\subseteq \rho_{\alpha}^{-1}(g)$. On the other hand, for the $C^*$-system $(C_{r}^{*}(\mathcal{R}),G,\alpha)$, there is a canonical dual coaction $\widehat{\alpha}: C_{r}^{*}(\mathcal{R})\rtimes_{\alpha,r} G\rightarrow (C_{r}^{*}(\mathcal{R})\rtimes_{\alpha,r} G)\otimes C^*(G)$ of $G$ on $C_{r}^{*}(\mathcal{R})\rtimes_{\alpha,r} G$, defined by $\widehat{\alpha}(a)=a\otimes I$ and $\widehat{\alpha}(u_g)=u_g\otimes v_g$, where $\{a: a\in C_{r}^{*}(\mathcal{R})\} \cup \{u_g:g\in G\}$ is the canonical generators of $C_{r}^{*}(\mathcal{R})\rtimes_{\alpha,r} G$, and   $C^*(G)$ is the full group $C^*$-algebra with generators $\{v_g: g\in G\}$ (\cite{Kaliszewski2016}). We conjecture that two systems $(C_r^*(\mathcal{R}\rtimes_{\alpha} G), G, \delta_{\alpha})$ and $(C_{r}^{*}(\mathcal{R})\rtimes_{\alpha,r} G, G,\widehat{\alpha})$ are conjugate when $G$ is amenable.
\end{remark}

\section{Continuous orbit equivalence of automorphism systems}

Given an automorphism system $G\curvearrowright_{\alpha} (X,\mathcal{R})$ on a compact metrizable space $X$, for $x\in X$, we let $[x]_G\, :=\{gx:\, g\in G\}$ and $[x]_{\mathcal{R}}\, :=\{y\in X:\, (x,y)\in \mathcal{R}\}$
be the orbits of $x$ under the action $\alpha$ and the relation $\mathcal{R}$, respectively. We call the set $[x]_{G,\mathcal{R}}=\{y\in X:\, (gx,y)\in \mathcal{R} \hbox{ for some $g\in G$}\}$ \emph{the bi-orbit} of $x$. Clearly, $[x]_{G,\mathcal{R}}=\cup_{y\in [x]_G}[y]_{\mathcal{R}}=\cup_{y\in [x]_{\mathcal{R}}}[y]_G=d((\mathcal{R}\rtimes_{\alpha} G)^x)=r((\mathcal{R}\rtimes_{\alpha} G)_x)$.

Recall that $G\curvearrowright X$ and $H\curvearrowright Y$ are {\sl orbit equivalent} if there  exists a homeomorphism $\varphi: X \rightarrow Y$ such that $\varphi([x]_G)=[\varphi(x)]_H$ for $x\in X$. They are said to be {\sl continuously orbit equivalent} if there exist a homeomorphism $\varphi: X \rightarrow Y$ and continuous maps $a: G\times X \rightarrow H$ and $b: H\times Y \rightarrow G$	such that $\varphi(gx)=a(g,x)\varphi(x)$ for $x\in X$ and $g\in G$, and
$\varphi^{-1}(hy)=b(h,y)\varphi^{-1}(y)$ for $y\in Y$ and $h\in H$ (\cite{Li2018}). Motivated by these notions, we introduce the following definitions.

\begin{definition}\quad
	Two systems $G\curvearrowright (X,\mathcal{R})$ and $H\curvearrowright (Y,\mathcal{S})$ are \emph{orbit equivalent} if there exists a homeomorphism $\varphi:\, X\rightarrow Y$ such that $\varphi([x]_{G, \mathcal{R}})=[\varphi(x)]_{H,\mathcal{S}}$ for $x\in X$.
	
\end{definition}

In this case, for $x,y\in X$ and $g\in G$ with $(gx,y)\in \mathcal{R}$, there exists $h$ in $H$ such that $(h\varphi(x), \varphi(y))\in \mathcal{S}$. Similarly, for $u,v\in Y$ and $h\in H$ with $(hu,v)\in \mathcal{S}$, there exists $g$ in $G$ such that $(g\varphi^{-1}(u), \varphi^{-1}(v))\in \mathcal{R}$. Thus, we have the following notion.

\begin{definition}\quad
	Two systems $G\curvearrowright (X,\mathcal{R})$ and $H\curvearrowright (Y,\mathcal{S})$ are \emph{continuously orbit equivalent}, we write $G\curvearrowright (X,\mathcal{R})\,\sim_{coe}\, H\curvearrowright (Y,\mathcal{S})$, if there exist a homeomorphism $\varphi:\, X\rightarrow Y$ and continuous maps $a:\,\, \mathcal{R}\times G\rightarrow H$ and  $b:\,\, \mathcal{S}\times H\rightarrow G$ such that
	the following maps:
	$$((x,y),g)\in \mathcal{R}\times G\rightarrow (\varphi(x), a((x,y),g)\varphi(g^{-1}y))\in \mathcal{S} $$
	and
	$$((x,y),g)\in \mathcal{S}\times H\rightarrow (\varphi^{-1}(x), b((x,y),g)\varphi^{-1}(g^{-1}y))\in \mathcal{R} $$
	are continuous.
\end{definition}

Clearly, continuous orbit equivalence implies orbit equivalence for automorphism systems. Assume a system $G\curvearrowright_{\alpha} X$ is free in the sense that, for $g\in G$ and $x\in X$, $gx=x$ only if $g=e$. We consider two automorphism systems $G\curvearrowright (X,\mathcal{R}_1)$ and $G\curvearrowright (X,\mathcal{R}_2)$, where $\mathcal{R}_1=\{(x,x):\, x\in X\}$ is the trivial \'{e}tale equivalence relation on $X$ under the relative product topology and $\mathcal{R}_2=\{(x,gx):\, x\in X, g\in G\}$ is the orbit equivalence relation under $\alpha$. Noticing that the map $(x,g)\in X\rtimes G\rightarrow (x, g^{-1}x)\in \mathcal{R}_2$ is a bijection, we transfer the product topology on $X\rtimes G$ over $\mathcal{R}_2$ via this map. Then $\mathcal{R}_2$ is an \'{e}tale equivalence relation on $X$.

\begin{proposition}\quad Assume that $G\curvearrowright X$ is free. Then $G\curvearrowright (X,\mathcal{R}_1)$ and $G\curvearrowright (X,\mathcal{R}_2)$ are  continuously orbit equivalent, but not conjugate.
\end{proposition}

\begin{proof}\quad Let $\varphi$ be the identity map on $X$, and let $a((x,x),g)=g$ for $((x,x),g)\in \mathcal{R}_1\times G$. For each $(x,y)\in \mathcal{R}_2$, there exists unique an element in $G$, denoted by $k(x,y)$,  such that $y=k(x,y)x$. Let $b((x,y),g)=k(x,y)^{-1}g$ for $((x,y),g)\in \mathcal{R}_2\times G$. Then $\varphi$, $a$ and $b$ satisfy the requirements in Definition 4.2, thus $G\curvearrowright (X,\mathcal{R}_1)$ and $G\curvearrowright (X,\mathcal{R}_2)$ are  continuously orbit equivalent.
	
	Since $\mathcal{R}_1$ and $\mathcal{R}_2$ are never isomorphic, $G\curvearrowright (X,\mathcal{R}_1)$ and $G\curvearrowright (X,\mathcal{R}_2)$ are not conjugate.

\end{proof}

Using the semi-direct product groupoid $\mathcal{R}\rtimes_{\alpha} G$ and the canonical homeomorphism $\gamma_0$, one can check the following lemma.

\begin{lemma}\quad
	Two systems $G\curvearrowright (X,\mathcal{R})\,\sim_{coe}\, H\curvearrowright (Y,\mathcal{S})$ if and only if there exist a homeomorphism $\varphi:\, X\rightarrow Y$ and continuous maps $a:\,\, \mathcal{R}\rtimes_{\alpha} G\rightarrow H$ and  $b:\,\, \mathcal{S}\rtimes_{\beta} H\rightarrow G$ such that
	the following maps:
	$$\Psi:\,\,(x,g,y)\in \mathcal{R}\rtimes_{\alpha} G\rightarrow (\varphi(x), a(x,g,y), \varphi(y))\in \mathcal{S}\rtimes_{\beta} H \eqno{(4.1)}$$
	and
	$$\widetilde{\Psi}:\, (u,h,v)\in \mathcal{S}\rtimes_{\beta} H\rightarrow (\varphi^{-1}(u), b(u,h,v),\varphi^{-1}(v))\in \mathcal{R}\rtimes_{\alpha} G \eqno{(4.2)}$$
	are continuous.
\end{lemma}

Recall that an \'{e}tale groupoid $\mathcal{G}$ is \emph{topologically principal} if $\{u\in \mathcal{G}^{(0)}:\,\, \mathcal{G}^u_u=\{u\}\}$ is dense in $\mathcal{G}^{(0)}$.
Since $\mathcal{G}$ is assumed to be second countable, it follows from \cite{Brown2014,Renault2008} that it is topologically principal if and only if the interior of $\mathcal{G}'$ is $\mathcal{G}^{(0)}$, where $\mathcal{G}'=\cup_{u\in \mathcal{G}^{(0)}}\mathcal{G}_u^u$ is the isotropy bundle of $\mathcal{G}$. For $G\curvearrowright_{\alpha} (X,\mathcal{R})$, we have
$$(\mathcal{R}\rtimes_{\alpha}G)'=\{(x,g,x):\,\, x\in X, g\in G, \,(x,gx)\in \mathcal{R}\}$$
and $$(\mathcal{R}\times_{\alpha}G)'=\{((x,gx),g):\,\, x\in X, g\in G, \,(x,gx)\in \mathcal{R}\}.$$
Moreover, we have that $\gamma_0((\mathcal{R}\rtimes_{\alpha}G)')=(\mathcal{R}\times_{\alpha}G)'$. Motivated by \cite{Hou2021,Matsumoto2019}, we have the following notion.

\begin{definition}\quad A system $G\curvearrowright (X,\mathcal{R})$ is said to be \emph{essentially free} if for every $e\neq g\in G$, $\{x\in X:\,\,  (x,gx)\notin \mathcal{R}\}$ is dense in $X$.
\end{definition}

One can easily see that  $G\curvearrowright (X,\mathcal{R})$ is essentially free, if and only if the interior of $\{x\in X:\,\, g[x]_{\mathcal{R}}=[x]_{\mathcal{R}}\}$ in $X$ is empty for every $g\neq e$.

\begin{lemma}\quad
	A system $G\curvearrowright_{\alpha} (X,\mathcal{R})$ is essentially free if and only if $\mathcal{R}\times_{\alpha} G$ ( or $\mathcal{R}\rtimes_{\alpha} G$) is topologically principal.
	
	Moreover, one of these conditions implies that both systems $G\curvearrowright X$ and $G\curvearrowright \mathcal{R}$ are topologically free.

\end{lemma}

\begin{proof}\quad It follows from the definitions that the topological principality of $\mathcal{R}\times_{\alpha} G$ implies the essential freeness of $G\curvearrowright (X,\mathcal{R})$, thus implies the topological freeness of $G\curvearrowright X$. To see that the essential freeness of $G\curvearrowright_{\alpha} (X,\mathcal{R})$ implies the topological principality of $\mathcal{R}\times_{\alpha} G$, we only need to show that $((x,gx),g)$ is not in the interior of $(\mathcal{R}\times_{\alpha}G)'$ in $(\mathcal{R}\times_{\alpha}G)$ for each $e\neq g\in G$ and $x\in X$ with $(x,gx)\in \mathcal{R}$.
	
	In fact, for otherwise, choose $e\neq g_0\in G$ and $x_0\in X$ such that $(x_0,g_0x_0)\in \mathcal{R}$ and $((x_0,g_0x_0),g_0)$ is an interior point of $(\mathcal{R}\times_{\alpha} G)'$. Then there exists an open neighbourhood  $\widetilde{U}$ of $(x_0,g_0x_0)$ in $\mathcal{R}$ such that $$((x_0,g_0x_0),g_0)\in \widetilde{U}\times \{g_0\}\subseteq (\mathcal{R}\times_{\alpha} G)'.$$ The last inclusion implies that $y=g_0x$ for each $(x,y)\in \widetilde{U}$. Hence $\{x\in X:\, (x,g_0x)\in \mathcal{R}\}$  contains the non-empty open subset $r(\widetilde{U})$ of $X$, which is contrast to the essential freeness of $G\curvearrowright (X,\mathcal{R})$.
	
	Assume $G\curvearrowright (X,\mathcal{R})$ is essentially free. Given $e\neq g\in G$ and a  non-empty open subset $U\subseteq \mathcal{R}$, it follows from the openness of $r(U)$ that there exists $x_0\in r(U)$ with $(x_0,gx_0)\notin \mathcal{R}$, thus $x_0\neq gx_0$. Choose $(x_0,y_0)\in U$. Then $g(x_0,y_0)\neq (x_0,y_0)$, which implies that $\{(x,y)\in \mathcal{R}:\, g(x,y)\neq (x,y)\}$ is dense in $\mathcal{R}$. Hence $G\curvearrowright \mathcal{R}$ is topologically free.
	
\end{proof}

\begin{remark}\quad The topological freeness of neither $G\curvearrowright X$ nor $G\curvearrowright \mathcal{R}$ can imply the essential freeness of $G\curvearrowright (X,\mathcal{R})$. To see this, if $G\curvearrowright X$ is free, then both systems $G\curvearrowright \mathcal{R}_1$ and $G\curvearrowright \mathcal{R}_2$ in Proposition 4.3  are free, and $G\curvearrowright (X,\mathcal{R}_1)$ is essentially free, but $G\curvearrowright (X,\mathcal{R}_2)$ is not.

	If $G\curvearrowright_{\alpha} (X,\mathcal{R})$ and $H\curvearrowright_{\beta} (Y,\mathcal{S})$ are essentially free, then the mappings $a$ and $b$ in Lemma 4.4 (or in Definition 4.2) are uniquely determined by (4.1) and (4.2). In fact, suppose that $a':\,\, \mathcal{R}\rtimes_{\alpha} G\rightarrow H$ is another continuous map such that
	$\Psi':\,\,(x,g,y)\in \mathcal{R}\rtimes_{\alpha} G\rightarrow (\varphi(x), a'(x,g,y),\varphi(y))\in \mathcal{S}\rtimes_{\beta} H$ is continuous. Then $$(x,g,y)\in \mathcal{R}\rtimes_{\alpha} G\rightarrow (a(x,g,y)\varphi(y), a'(x,g,y)\varphi(y))\in \mathcal{S}$$ is continuous. Hence, from the continuity of $a$, $a'$ and $\rho_{\alpha}$, for $(x,g,y)\in \mathcal{R}\rtimes_{\alpha} G$, there exists an open neighbourhood $\widetilde{U}$ of $(x,g,y)$ such that the map $d|_{\widetilde{U}}:\, \widetilde{U}\rightarrow d(\widetilde{U})$ is a homeomorphism, $\rho_{\alpha}(\gamma)=g$,
	$a(\gamma)=a(x,g,y)$, and $a'(\gamma)=a'(x,g,y)$ for each $\gamma\in \widetilde{U}$.
	For each $z\in \varphi(d(\widetilde{U}))$, choose $\gamma\in \widetilde{U}$ such that $z=\varphi(d(\gamma))$. The choice of  $\widetilde{U}$ implies that we can assume that $\gamma=(u,g,v)$, thus $z=\varphi(v)$. Note that $(\varphi(u),a(\gamma)z)$ and $(\varphi(u),a'(\gamma)z)$, thus $(a(\gamma)z, a'(\gamma)z)$ are in $\mathcal{S}$. Hence $(a(x,g,y)z, a'(x,g,y)z)\in \mathcal{S}$ for each $z\in \varphi(d(\widetilde{U}))$. The essential freeness of $H\curvearrowright_{\beta} (Y,\mathcal{S})$ implies $a(x,g,y)=a'(x,g,y)$. By symmetry, $b$ is uniquely determined by (4.2).
\end{remark}

\begin{lemma}\quad In Definition 4.2, if $G\curvearrowright_{\alpha} (X,\mathcal{R})$ and $H\curvearrowright_{\beta} (Y,\mathcal{S})$ are essentially free, then the mappings $a$ and $b$ are cocycles on $\mathcal{R}\times_{\alpha} G$ and $\mathcal{S}\times_{\beta} H$, respectively.
\end{lemma}

\begin{proof}\quad We only need to show that the mappings $a$ and $b$ in Lemma 4.4 are cocycles.
	Let $\gamma_1=(x,g,y), \gamma_2=(y,h,z)\in \mathcal{R}\rtimes_{\alpha}G$ be arbitrary, and write $\gamma'=\gamma_1\gamma_2=(x,gh,z)$. From the continuity of $a$ and $\rho_{\alpha}$, choose open neighbourhoods $U$, $V$ and $W$ of $\gamma_1$, $\gamma_2$ and $\gamma'$ in $\mathcal{R}\rtimes_{\alpha}G$, respectively, such that
	$a(\gamma)=a(\gamma_1),\, \rho_{\alpha}(\gamma)=g$ for each $\gamma\in U$, $a(\eta)=a(\gamma_2),\, \rho_{\alpha}(\eta)=h$ for each $\eta\in V$, and $a(\sigma)=a(\gamma'),\, \rho_{\alpha}(\sigma)=gh$ for each $\sigma\in W$. Since the multiplication on $(\mathcal{R}\rtimes_{\alpha}G)^{(2)}$ is continuous at $(\gamma_1,\gamma_2)$, we can  assume that $\gamma\eta\in W$ when $\gamma\in U$, $\eta\in V$ and $(\gamma,\eta)\in (\mathcal{R}\rtimes_{\alpha}G)^{(2)}$. Also since the range $r$ and domain $d$ are local homeomorphisms and  $d(\gamma_1)=r(\gamma_2)=y$, we can assume that the restrictions $d|_U$ and $r|_V$ are homeomorphisms onto their respective ranges and $d(U)=r(V)$.
	
	For each $\widetilde{y}\in \varphi(d(V))$, choose $\eta\in V$ such that $\widetilde{y}=\varphi(d(\eta))$. The choice of $V$ permits us to assume that $\eta=(v,h,w)$ and $a(\eta)=a(\gamma_2)$. Hence $\widetilde{y}=\varphi(w)$. Since $v\in r(V)=d(U)$, it follows from the choice of $U$ that we have a $\gamma=(u,g,v)\in U$ and $a(\gamma)=a(\gamma_1)$. Hence $\gamma\eta=(u,gh,w)\in W$ and $a(\gamma\eta)=a(\gamma')$.
	The hypothesis on $\Psi$ in Lemma 4.4 implies that $(\varphi(u), a(\gamma)\varphi(v))$, $(\varphi(v), a(\eta)\varphi(w))$ and $(\varphi(u), a(\gamma\eta)\varphi(w))$ are all  in $\mathcal{S}$. Thus, $(a(x,g,y)a(y,h,z)\widetilde{y}, a(x,gh,z)\widetilde{y})$ is in $\mathcal{S}$ for every $\widetilde{y}\in \varphi(d(V))$. The essential freeness of $H\curvearrowright_{\beta}(Y,\mathcal{S})$ implies $a(x,g,y)a(y,h,z)=a(x,gh,z)$, thus $a$ is a cocycle. By a similar way, we can show that $b$ is a cocycle.
\end{proof}

\begin{lemma}\quad In Definition 4.2, if $G\curvearrowright_{\alpha} (X,\mathcal{R})$ and $H\curvearrowright_{\beta} (Y,\mathcal{S})$ are essentially free, then
	$$b((\varphi(x),a((x,y),g)\varphi(g^{-1}y)),a((x,y),g))=g,$$
	$$a((\varphi^{-1}(u),b((u,v),h)\varphi^{-1}(h^{-1}v)),b((u,v),h))=h$$
	for every $((x,y),g)\in \mathcal{R}\times G$ and $((u,v),h)\in \mathcal{S}\times H$.
\end{lemma}

\begin{proof}\quad We only show that the maps $a$ and $b$ in Lemma 4.4 satisfy that
	$$b(\varphi(x),a(x,g,y),\varphi(y))=g,\,\,\, a(\varphi^{-1}(u),b(u,h,v),\varphi^{-1}(v))=h$$
	for every $(x,g,y)\in \mathcal{R}\rtimes_{\alpha} G$ and $(u,h,v)\in \mathcal{S}\rtimes_{\beta} H$.
	
	As before, let $\rho_{\alpha}$ and $\rho_{\beta}$ be the canonical cocycles on $\mathcal{R}\rtimes_{\alpha} G$ and $\mathcal{S}\rtimes_{\beta} H$, respectively. For an arbitrary $(x,g,y)\in \mathcal{R}\rtimes_{\alpha} G$, we have
	$(\varphi(x),h, \varphi(y))\in \mathcal{S}\rtimes_{\beta} H$, where $h=a(x,g,y)$.
	From the continuity of $b$ and $\rho_{\beta}$, there exists an open neighbourhood $U$ of $(\varphi(x),h, \varphi(y))$ in $\mathcal{S}\rtimes_{\beta} H$ such that $\rho_{\beta}(\gamma)=h$, $b(\gamma)=b(\varphi(x),h,\varphi(y))$ for every $\gamma\in U$, and $r|_{U}$, $d|_U$ are homeomorphisms from $U$ onto $r(U)$ and $d(U)$, respectively.
	
	By the continuity of $\rho_{\alpha}$, $\Psi$ and $a$ at $(x,g,y)$, as well as that of $\varphi$ at $x$ and $y$, there is an open neighbourhood $V$ of  $(x,g,y)$ in $\mathcal{R}\rtimes_{\alpha} G$ such that
	\begin{enumerate}
		\item[(i)] $\rho_{\alpha}(\gamma)=g$, $a(\gamma)=h$ and $\Psi(\gamma)\in U$ for every $\gamma\in V$;
		\item[(ii)] $r|_V$ and $d|_V$ are homeomorphisms from $V$ onto $r(V)$ and $d(V)$, respectively;
		\item[(iii)] $\varphi(r(V))\subseteq r(U)$ and $\varphi(d(V))\subseteq d(U)$.
	\end{enumerate}
	
	For each $v\in d(V)$, let $\gamma\in V$ such that $d(\gamma)=v$. The above condition (i) implies that we can let $\gamma=(u,g,v)$ and have $a(\gamma)=h$, thus
	$\Psi(\gamma)=(\varphi(u),h, \varphi(v))\in U$. The map $\widetilde{\Psi}$ gives that $(u, b(\varphi(u),h, \varphi(v))v)\in \mathcal{R}$. From the choice of $U$, $b(\varphi(u),h,\varphi(v))=b(\varphi(x),h,\varphi(y))$. So $(u, b(\varphi(x),h,\varphi(y))v)\in \mathcal{R}$.
	Also since $(u,gv)\in \mathcal{R}$, we have $(b(\varphi(x),h,\varphi(y))v,gv)\in \mathcal{R}$. The essential freeness of $G\curvearrowright_{\alpha} (X,\mathcal{R})$ implies that $b(\varphi(x),a(x,g,y),\varphi(y))=b(\varphi(x),h,\varphi(y))=g$.
	
	By a similar way, we can show that $a(\varphi^{-1}(u),b(u,h,v),\varphi^{-1}(v))=h$ for each $(u,h,v)\in \mathcal{S}\rtimes_{\beta} H$.
	
\end{proof}

The following definition comes from \cite[Definition 4.1]{Hou2021}.

\begin{definition}\quad For two \'{e}tale equivalence relations $\mathcal{R}$ and $\mathcal{S}$ on $X$ and $Y$, let $G\curvearrowright X$ and $H\curvearrowright Y$ be two systems generating two automorphism systems $G\curvearrowright (X,\mathcal{R})$ and $H\curvearrowright (Y,\mathcal{S})$. We say that $G\curvearrowright X$ and $H\curvearrowright Y$ \emph{continuously orbit equivalent up to $\mathcal{R}$ and $\mathcal{S}$}, if there exist a homoeomorphism, $\varphi: X\rightarrow Y$, continuous cocycles $a: X\rtimes G\rightarrow H$,
	$b: Y\rtimes H\rightarrow G$, $\sigma:\,\mathcal{R}\rightarrow H$, and $\tau:\, \mathcal{S}\rightarrow G$ satisfying the following conditions:
	\begin{enumerate}
		\item[(i)] $\sigma(x,y)a(y,g)=a(x,g)\sigma(g^{-1}x,g^{-1}y)$ for $(x,y)\in \mathcal{R}$ and $g\in G$;
		\item[(ii)] $\tau(x,y)b(y,g)=b(x,g)\tau(g^{-1}x,g^{-1}y)$ for $(x,y)\in \mathcal{S}$ and $g\in H$;
		\item[(iii)] The map, $\xi_1: (x,g)\in X\times G\rightarrow (a(x,g)^{-1}\varphi(x),\varphi (g^{-1}x))\in \mathcal{S}$, is well-defined and continuous; Moreover, $$b(\varphi(x),a(x,g))\,\tau(\xi_1(x,g))=g\,\,\mbox{ for $x\in X$ and $g\in G$}.$$
		\item[(iv)] The map, $\xi_2: (x,g)\in Y\times H\rightarrow (b(x,g)^{-1}\varphi^{-1} (x),\varphi^{-1} (g^{-1}x))\in \mathcal{R}$, is well-defined and continuous;  Moreover, $$a(\varphi^{-1}(x),b(x,g))\sigma(\xi_2(x,g))=g\, \mbox{ for $x\in Y$ and $g\in H$}.$$
		\item[(v)] The map, $\eta_1: (x,y)\in \mathcal{R}\rightarrow (\sigma(x,y)^{-1}\varphi(x),\varphi (y))\in \mathcal{S}$, is well-defined and continuous; Moreover, $$b(\varphi(x),\sigma(x,y))\,\tau(\eta_1(x,y))=e\, \mbox{ for $(x,y)\in \mathcal{R}$}.$$
		\item[(vi)] The map, $\eta_2: (x,y)\in \mathcal{S}\rightarrow (\tau(x,y)^{-1}\varphi^{-1}(x),\varphi^{-1} (y))\in \mathcal{R}$, is well-defined and continuous; Moreover, $$a(\varphi^{-1}(x),\tau(x,y))\,\sigma(\eta_2(x,y))=e\,\mbox{ for $(x,y)\in \mathcal{S}$}.$$
	\end{enumerate}
\end{definition}

\begin{proposition}\quad Let $G\curvearrowright_{\alpha} (X,\mathcal{R})$ and $H\curvearrowright_{\beta} (Y,\mathcal{S})$ be two automorphism systems. Then $G\curvearrowright X$ and $H\curvearrowright Y$ are continuously orbit equivalent up to $\mathcal{R}$ and $\mathcal{S}$ if and only if $\mathcal{R}\rtimes_{\alpha} G$ and $\mathcal{S}\rtimes_{\beta} H$ are isomorphic as \'{e}tale groupoids.
\end{proposition}

The proof of this proposition is the same as that of \cite[Theorem 4.2]{Hou2021} in which the local conjugacy is not necessary. We only provide a brief proof. For details, see \cite[Theorem 4.2]{Hou2021}.

\begin{proof}\quad
	Assume that $\Lambda:\,\,\mathcal{R}\rtimes_{\alpha} G\rightarrow \mathcal{S}\rtimes_{\beta} H$ is an isomorphism. Let $\varphi$ be the restriction of $\Lambda$ to the unit space $(\mathcal{R}\rtimes G)^{(0)}$ and let $a(x,g)=\rho_{\beta}\Lambda(x,g,g^{-1}x)$, $\sigma(x,y)=\rho_{\beta}\Lambda(x,e,y)$, and $b(u,h)=\rho_{\alpha}\Lambda^{-1}(u,h,h^{-1}u)$, $\tau(u,v)=\rho_{\alpha}\Lambda^{-1}(u,e,v)$. Then $\varphi$, $a,b$, $\sigma$ and $\tau$ satisfy the requirements in Definition 4.10, thus $G\curvearrowright X$ and $H\curvearrowright Y$ are continuously orbit equivalent up to $\mathcal{R}$ and $\mathcal{S}$.
	
	Conversely, assume that there are maps $\varphi$, $a,b$, $\sigma$ and $\tau$ satisfying the requirements in Definition 4.10. Define
	$$\Lambda(x,g,y)=(\varphi(x), a(x,g)\sigma(g^{-1}x,y),\varphi(y))\,\,\mbox{ for $(x,g,y)\in \mathcal{R}_{\alpha}\rtimes G$.}$$
	Then $\Lambda$ is an isomorphism from $\mathcal{R}\rtimes_{\alpha} G$ onto $\mathcal{S}\rtimes_{\beta} H$, whose inverse $\Lambda^{-1}$ is defined by $\Lambda^{-1}(u,h,v)=(\varphi^{-1}(u), b(u,h)\tau(h^{-1}u,v),\varphi^{-1}(v))$.
\end{proof}

\begin{theorem}\quad Assume that $G\curvearrowright_{\alpha} (X,\mathcal{R})$ and $H\curvearrowright_{\beta} (Y,\mathcal{S})$ are essentially free. Then the following statements are equivalent.
	\begin{enumerate}
		\item[(i)] $G\curvearrowright_{\alpha} (X,\mathcal{R}) \sim_{coe} H\curvearrowright_{\beta} (Y,\mathcal{S})$;
		\item[(ii)] $G\curvearrowright_{\alpha} X$ and $H\curvearrowright_{\beta} Y$ are continuously orbit equivalent up to $\mathcal{R}$ and $\mathcal{S}$;
		\item[(iii)] $\mathcal{R}\rtimes_{\alpha} G$ and $\mathcal{S}\rtimes_{\beta} H$ are isomorphic as \'{e}tale groupoids;
		\item[(iv)] there exists a $C^*$-isomorphism $\Phi$ from $C_r^*(\mathcal{R}\rtimes_{\alpha} G)$ onto $C_r^*(\mathcal{S}\rtimes_{\beta} H)$ such that $\Phi(C(X))=C(Y)$.
	\end{enumerate}
\end{theorem}

\begin{proof}\quad The equivalence of (ii) and (iii) follows from Proposition 4.11. From Lemma 4.6, $\mathcal{R}\rtimes_{\alpha} G$ and $\mathcal{S}\rtimes_{\beta} H$ are topological principal, thus the equivalence of (iii) and (iv) follows from \cite{Carlsen2021,Renault2008}.
	
	Assume (iii) holds, i.e., there is an isomorphism $\Lambda$ from $\mathcal{R}\rtimes_{\alpha} G$ onto $\mathcal{S}\rtimes_{\beta} H$. Let $\varphi$ be the restriction of $\Lambda$ to the unit space $X$, and let $a(x,g,y)=\rho_{\beta}\Lambda(x,g,y)$ for $(x,g,y)\in \mathcal{R}\rtimes_{\alpha} G$, $b(u,h,v)=\rho_{\alpha}\Lambda^{-1}(u,h,v)$ for $(u,h,v)\in \mathcal{S}\rtimes_{\beta} H$. Then $\varphi$ is a homeomorphism from $X$ onto $Y$, and $\Lambda(x,g,y)=(\varphi(x),a(x,g,y),\varphi(y))$ and $\Lambda^{-1}(u,h,v)=(\varphi^{-1}(u),b(u,h,v),\varphi^{-1}(v))$.
	So $\varphi$, $a$ and $b$ satisfy the requirements in Lemma 4.4, thus $G\curvearrowright_{\alpha} (X,\mathcal{R}) \sim_{coe} H\curvearrowright_{\beta} (Y,\mathcal{S})$, i.e., (i) holds.
	
	Assume (i) holds. From Lemma 4.4, there exist mappings $\varphi$, $a$ and $b$ such that the mappings
	$\Psi:\, (x,g,y)\in \mathcal{R}\rtimes_{\alpha} G\rightarrow (\varphi(x), a(x,g,y),\varphi(y))\in \mathcal{S}\rtimes_{\beta} H$ and $\widetilde{\Psi}: (u,h,v)\in \mathcal{S}\rtimes_{\beta} H\rightarrow (\varphi^{-1}(u), b(u,h,v),\varphi^{-1}(v))\in \mathcal{R}\rtimes_{\alpha} G$ are continuous. From Lemma 4.8 and Lemma 4.9, $\Psi$ and $\widetilde{\Psi}$ are \'{e}tale groupoid isomorphisms and inverse to each other, thus (iii) holds.

\end{proof}

\begin{remark}\quad For $G\curvearrowright_{\alpha} X$, let $\mathcal{R}_1=\{(x,x):\, x\in X\}$  be as in Proposition 4.3. Then $\mathcal{R}_1\rtimes_{\alpha} G$ is isomorphic to the transformation groupoid $X\rtimes G$, and the notions of continuous orbit equivalence for $G\curvearrowright (X,\mathcal{R}_1)$ and $G\curvearrowright X$ in the Li's sense are consistent. Hence Theorem 4.12 is a generalization of Theorem 1.2 in \cite{Li2018}.

\end{remark}

There are two special cases for orbit equivalence of two systems $G\curvearrowright_{\alpha} (X,\mathcal{R})$ and $H\curvearrowright_{\beta} (Y,\mathcal{S})$ via a homeomorphism $\varphi:\, X\rightarrow Y$. One is, for each $g\in G$, there is $h\in H$ such that $(h\varphi(x),\varphi(y))\in \mathcal{S}$ for each $x,y\in X$ with $(gx,y)\in \mathcal{R}$, and by symmetry, for each $h\in H$, there is $g\in G$ such that $(g\varphi^{-1}(x),\varphi^{-1}(y))\in \mathcal{R}$ for each $x,y\in Y$ with $(hx,y)\in \mathcal{S}$. The other is, for each $g\in G$ and $x\in X$, there is $h\in H$ such that $(h\varphi(x),\varphi(z))\in \mathcal{S}$ for each $(gx,z)\in \mathcal{R}$, and by symmetry, for each $h\in H$ and $y\in Y$, there is $g\in G$ such that $(g\varphi^{-1}(y),\varphi^{-1}(z))\in \mathcal{R}$ for each $(hy,z)\in \mathcal{S}$. Inspired by these ideas, we have the following notions, comparing with those of (strong) asymptotic conjugation in \cite[Definition 4.4]{Hou2021}.

\begin{definition}\quad
	We say $G\curvearrowright (X,\mathcal{R})$ and $H\curvearrowright (Y,\mathcal{S})$ \emph{strongly continuously orbit equivalent}, write $G\curvearrowright (X,\mathcal{R})\sim_{scoe} H\curvearrowright (Y,\mathcal{S})$, if they are continuously orbit equivalent and in Definition 4.2 we can take the maps $a(\gamma,g)=a(\gamma',g)$ and $b(\nu,h)=b(\nu',h)$ for all $\gamma,\gamma'\in \mathcal{R}$ and $\nu,\nu'\in \mathcal{S}$.

	We say these two systems \emph{weakly continuously orbit equivalent}, write $G\curvearrowright (X,\mathcal{R})\sim_{wcoe} H\curvearrowright (Y,\mathcal{S})$, if they are continuously orbit equivalent and in Definition 4.2 we can take the maps $a(\gamma,g)=a(\gamma',g)$ for $\gamma,\gamma'\in \mathcal{R}$ with $d(\gamma)=d(\gamma')$, and $b(\nu,h)=b(\nu',h)$ for $\nu,\nu'\in \mathcal{S}$ with $d(\nu)=d(\nu')$.

\end{definition}

\begin{remark}\quad
	Clearly, the strong continuous orbit equivalence implies the weak one.  If $G\curvearrowright_{\alpha} X$ is free, then $G\curvearrowright (X,\mathcal{R}_1)$ and $G\curvearrowright (X,\mathcal{R}_2)$ in Proposition 4.3 are continuously orbit equivalent, but not weakly continuously orbit equivalent, because they do  not satisfy the second special case.
\end{remark}

The following corollary is an analogy to \cite[Proposition 4.5]{Hou2021}

\begin{corollary}\quad Assume that $G\curvearrowright_{\alpha} (X,\mathcal{R})$ and $H\curvearrowright_{\beta} (Y,\mathcal{S})$ are essentially free. Then
	\begin{enumerate}
		\item[(i)] $G\curvearrowright_{\alpha} (X,\mathcal{R}) \sim_{wcoe} H\curvearrowright_{\beta} (Y,\mathcal{S})$ if and only if there is an isomorphism $\Lambda:\,\,\mathcal{R}\rtimes_{\alpha} G\rightarrow \mathcal{S}\rtimes_{\beta} H$ such that $\Lambda(\mathcal{R})=\mathcal{S}$.  
		
		Moreover, if these conditions hold, then there is a $C^*$-isomorphism $\Phi: C_r^*(\mathcal{R}\rtimes_{\alpha} G)\rightarrow C_r^*(\mathcal{S}\rtimes_{\beta} H)$ such that $\Phi(C(X))=C(Y)$ and $\Phi(C_r^*(\mathcal{R}))=C_r^*(\mathcal{S})$.
		
		\item[(ii)] $G\curvearrowright_{\alpha} (X,\mathcal{R}) \sim_{scoe} H\curvearrowright_{\beta} (Y,\mathcal{S})$ if and only if   there exist a homeomorphism $\varphi:\, X\rightarrow Y$ and a group isomorphism $\theta: G\rightarrow H$ such that $\Lambda:\, (x,g,y)\in \mathcal{R}\rtimes_{\alpha} G\rightarrow (\varphi(x),\theta(g),\varphi(y))\in \mathcal{S}\rtimes_{\beta} H$ is an isomorphism if and only if   there exist an \'{e}tale groupoid isomorphism $\Lambda:\, \mathcal{R}\rtimes_{\alpha} G\rightarrow \mathcal{S}\rtimes_{\beta} H$ and a group isomorphism $\theta: G\rightarrow H$ such that $\theta\rho_{\alpha}=\rho_{\beta}\Lambda$ if and only if two coaction systems $(C_r^*(\mathcal{R}\rtimes_{\alpha} G), G, \delta_{\alpha})$ and $(C_r^*(\mathcal{S}\rtimes_{\alpha} H), H, \delta_{\beta})$ are conjugate by a conjuagcy $\phi$ with $\phi(C(X))=C(Y)$.

	\end{enumerate}
	
	Furthermore,  when $\mathcal{R}$ and $\mathcal{S}$ are minimal or $X$ and $Y$ are connected, these two notions of strong continuous orbit equivalence and weak continuous orbit equivalence are consistent.
\end{corollary}

\begin{proof}\quad
	One can check that if the map $a$ in Definition 4.2 is a cocycle on $\mathcal{R}\times_{\alpha} G$, then $a(\gamma,g)=a(\gamma',g)$ for $\gamma,\gamma'\in \mathcal{R}$ with $d(\gamma)=d(\gamma')$ if and only if $a(\gamma,e)=e$ for all $\gamma\in \mathcal{R}$. By symmetry, $b$ has a similar characterization when it is a cocycle. From Lemma 4.8, Theorem 4.12 and its proof, we can obtain (i) and (ii), where the equivalence of the last two statements in (ii) comes from \cite[Theorem 6.2]{Carlsen2021}.

	We now show that the weak continuous orbit equivalence of $G\curvearrowright_{\alpha} (X,\mathcal{R})$ and $H\curvearrowright_{\beta} (Y,\mathcal{S})$ implies the strong one when $\mathcal{R}$ and $\mathcal{S}$ are minimal or $X$ and $Y$ are connected. To see this, by assumption and the first paragraph of this proof, we have a homeomorphism $\varphi$ and two continuous cocycles $a, b$  with $a(x,e,y)=e$ for all $(x,y)\in \mathcal{R}$ and $b(u,e,v)=e$ for $(u,v)\in \mathcal{S}$, satisfying Lemma 4.4.
	
	Assume that $X$ and $Y$ are connected. For each $g\in G$, the map $x\in X \rightarrow a(x,g,g^{-1}x)\in H$ is continuous, thus it is a constant. Hence $a(x,g,g^{-1}x)=a(y,g,g^{-1}y)$ for all $x,y\in X$. By symmetry, $b$ has a similar property.
	
	Assume that $\mathcal{R}$ and $\mathcal{S}$ are minimal. For $(x,y)\in \mathcal{R}$ and $g\in G$, since $$(x,g,g^{-1}x)(g^{-1}x,e,g^{-1}y)(g^{-1}y,g^{-1},y)=(x,e,y),$$ we have $a(x,g,g^{-1}x)=a(g^{-1}y,g^{-1},y)^{-1}=a(y,g,g^{-1}y)$. Given arbitrary  $x,y\in X$ and $g\in G$, we choose a sequence $\{x_n\}$ in $[x]_{\mathcal{R}}$ converging to $y$ in $X$. Thus $\{(x_n,g,g^{-1}x_n)\}$ converges to $(y,g,g^{-1}y)$ in $\mathcal{R}\rtimes_{\alpha}G$, thus the continuity of $a$ implies that $a(x,g,g^{-1}x)=a(y,g,g^{-1}y)$.
	
	Remark that $a(x,g,y)=a(x,g,g^{-1}x)a(g^{-1}x,e,y)=a(x,g,g^{-1}x)$ for $(x,g,y)\in \mathcal{R}\rtimes_{\alpha} G$. Consequently, if one of the above two assumptions holds, then $a(x,g,y)=a(u,g,v)$ for $(x,g,y),\, (u,g,v)\in \mathcal{R}\rtimes_{\alpha} G$. By a similar way, we can show that $b$ satisfies a similar requirement. Hence $G\curvearrowright_{\alpha} (X,\mathcal{R})$ and $H\curvearrowright_{\beta} (Y,\mathcal{S})$ are strongly continuously orbit equivalent.

\end{proof}

\section{ Local conjugacy relations from expansive systems}

The condition of essential freeness of automorphism systems in Theorem 4.12 and Corollary 4.16 is necessary. In this section, we give some examples satisfying the requirement. Recall that a system $G\curvearrowright_{\alpha} X$ is called \emph{expansive} if the action $\alpha$ is expansive, which means for a metric $d$ on $X$ compatible with the topology, there exists a constant $\delta>0$ such that, for $x,y\in X$,  if $d(gx,gy)<\delta$ for all $g\in G$ then $x=y$. For convenience, given a real-valued function $\psi$ on $G$, the notation $\lim\limits_{g\rightarrow\infty}\psi(g)=0$ means that, for any $\epsilon>0$, there exists a finite subset $F$ of $G$ such that $|\psi(g)|<\epsilon$ for all $g\notin F$.

A  triple $(U,V,\gamma)$, consisting of open subsets $U$, $V$ of $X$ and a homeomorphism $\gamma:\, U\rightarrow V$, is called  a local conjugacy, if $\lim\limits_{g\rightarrow \infty}\sup_{z\in U}d(gz,g\gamma (z))=0$.
Two points  $x$ and $y$ in $X$ are  said to be \emph{locally conjugate}, if there exists a local conjugacy $(U,V,\gamma)$ such that $x\in U$, $y\in V$ and $\gamma (x)=y$. Let
$$\mathcal{R}_{\alpha}=\{(x,y)\in X\times X:\,\,\mbox{$x$ and $y$ are locally conjugate}\}$$
be the local conjugacy relation on $X$. From \cite{Ruelle1988} (also see \cite{Thomsen2010}), $\mathcal{R}_{\alpha}$ is an \'{e}tale  equivalence relation on $X$ under the topology whose base consists of the sets of the form
$$\{(x,\gamma (x)):\, x\in U\},$$
where $(U,V,\gamma)$ is a local conjugacy. Moreover, $G\curvearrowright_{\alpha} X$ induces an automorphism system $G\curvearrowright_{\alpha}\mathcal{R}_{\alpha}$: $g(x,y)=(gx,gy)$ for $g\in G$ and $(x,y)\in\mathcal{R}_{\alpha}$. Thus we have an automorphism system $G\curvearrowright_{\alpha} (X, \mathcal{R}_{\alpha})$.

\begin{remark}\quad
	If two expansive systems $G\curvearrowright_{\alpha} X$ and $H\curvearrowright_{\beta} Y$ are conjugate by a homeomorphism $\varphi$ from $X$ onto $Y$ and a group isomorphism $\rho$ from $G$ onto $H$, then $(\varphi(U),\varphi(V),\varphi\gamma\varphi^{-1}|_{\varphi(U)})$ is a local conjugacy for each local conjugacy $(U,V,\gamma)$, thus $\varphi\times\varphi:\, (x,y)\in \mathcal{R}_{\alpha}\rightarrow (\varphi(x),\varphi(y))\in \mathcal{R}_{\beta}$ is an isomorphism. Hence $G\curvearrowright_{\alpha} (X, \mathcal{R}_{\alpha})$ and $H\curvearrowright_{\beta} (Y, \mathcal{R}_{\beta})$ are conjugate, thus two notions of conjugacy for $G\curvearrowright_{\alpha} (X, \mathcal{R}_{\alpha})$ and $G\curvearrowright_{\alpha} X$ are consistent.
\end{remark}

From \cite{Hou2021} and \cite{Matsumoto2019}, the automorphism systems of local conjugacy relations associated to a full shift $G\curvearrowright A^G$ over a finite set $A$ and an irreducible Smale space $(X,\psi)$ are essentially free. The following result generalizes Matsumoto's result in the  Smale space case to the  $\mathbb{Z}$-expansive system case.

\begin{theorem}\quad	Let $\mathcal{R}_{\varphi}$ be the local conjugacy relation from an expansive system $\mathbb{Z}\curvearrowright_{\alpha} X$ generated by a  homeomorphism $\varphi$ on $X$. Assume that $X$ is infinite and has no isolated points. Then $\mathbb{Z}\curvearrowright_{\alpha} (X,\mathcal{R}_{\varphi})$ is essentially free.
\end{theorem}

\begin{proof}\quad
	For an arbitrary integer $p\geq 1$, we first claim that the set $$X_p=\{x\in X:\, \lim\limits_{n\rightarrow \infty} \varphi^{pn}(x)\,\, \hbox{ and }\,\,\lim\limits_{n\rightarrow \infty}\varphi^{-pn}(x)\, \hbox{ exist}\}$$ is countable.
	
	In fact, when $p=1$, it follows from \cite[Theorem 2.2.22]{Aoki1994} that $X_1$ is countable. Moreover, the expansiveness of $\varphi^p$ implies that $X_p$ is also countable for every $p$. For completion, we provide a proof for the claim. Since $\varphi$ is expansive, it follows that the $p$-periodic point set $F_p(\varphi)=\{x\in X:\, \varphi^p(x)=x\}$ is finite, say $F_p(\varphi)=\{y_1, y_2,\cdots, y_k\}$. For each $x\in X_p$, let $\lim\limits_{n\rightarrow \infty} \varphi^{pn}(x)=y$ and $\lim\limits_{n\rightarrow \infty}\varphi^{-pn}(x)=z$. One can see that $y,z\in F_p(\varphi)$. For $1\leq i,j\leq k$, set $X_p(i,j)=\{x\in X:\,\, \lim\limits_{n\rightarrow \infty} \varphi^{pn}(x)=y_i, \,\,\lim\limits_{n\rightarrow \infty}\varphi^{-pn}(x)=y_j\}$. Given $x\in X_p(i,j)$, we have $\lim\limits_{n\rightarrow \infty} \varphi^{pn+r}(x)=\varphi^r(y_i)$ and $\lim\limits_{n\rightarrow \infty}\varphi^{-pn-r}(x)=\varphi^{-r}(y_j)$ for each $0\leq r\leq p-1$. Hence there exists an integer $N\geq 2$ such that $d(\varphi^n(x), \varphi^n(y_i))<\frac c2$ for all $n\geq N$ and $d(\varphi^{-n}(x), \varphi^{-n}(y_j))<\frac c2$ for all $n\leq -N$, where $c$ is an expansive constant for $\varphi$. Set
	$$X_{p,N}(i,j)=\left\{x\in X:\, \begin{array}{l}
		d(\varphi^n(x), \varphi^n(y_i))<\frac c2 \hbox{ for all $n\geq N$}\\
		d(\varphi^{-n}(x), \varphi^{-n}(y_j))<\frac c2 \hbox{  for all $n\leq -N$}\end{array}\right\}.$$
	Thus $X_p(i,j)=\cup_{N\geq 2} X_{p,N}(i,j)$ and $X=\cup_{1\leq i,j\leq k}X_p(i,j)$. To finish the claim, we show, for each $N\geq 2$ and $1\leq i,j\leq k$, the set $X_{p,N}(i,j)$ is finite.
	
	For otherwise, $X_{p,N}(i,j)$ is infinite for some $i,j, N$. Choose $\delta<\frac c2$ such that if $d(y,z)\leq \delta$ for $y,z\in X$ then $d(\varphi^l(y), \varphi^l(z))<\frac c2$ for each integer $l$ with $|l|\leq N-1$. Since $X_{p,N}(i,j)$ is infinite, there are two different $y,z$ in $X_{p,N}(i,j)$ such that $d(y,z)<\delta$. Thus $d(\varphi^l(y), \varphi^l(z))<c$ for every integer $l$, which implies that $y=z$ by expansiveness of $\varphi$ and is a contradiction. We have established the claim.
	
	For each $p\geq 1$, we next claim that if $x\in X$ with $(x,\varphi^{p}(x))\in \mathcal{R}_{\varphi}$, then $x\in X_p$.
	
	We use the method in \cite[Lemma 5.3]{Putnam1999} to complete the claim. Assume that $z$ is a limit point of $\left\lbrace\varphi^{pn}(x) |n\geq 1\right\rbrace$. Choose a subsequence $\left\lbrace m_{n}\right\rbrace $ of positive integers such that $\underset{n\rightarrow \infty}{\mathrm{lim}}d(z, \varphi^{pm_{n}}(x))=0$. Thus $\lim\limits_{n\rightarrow \infty}d(\varphi^{p}(z), \varphi^{pm_{n}}(\varphi^{p}(x)))=0$.
	Since $(x,\varphi^{p}(x))\in \mathcal{R}_{\varphi}$, we have  $\underset{|n|\rightarrow \infty}{\mathrm{lim}} d(\varphi^{n}(x),\varphi^{n}(\varphi^{p}(x)))=0$.  Consequently, $\varphi^{p}(z)=z$, which implies that
	each limit point of $\left\lbrace\varphi^{pn}(x) |n\geq1\right\rbrace$ is in $F_p(\varphi)=\{y_1, y_2,\cdots, y_k\}$.
	Choose an open neighbourhood $U_{i}$ of $y_{i}$ such that $\overline{U_{i}}\cap \overline{U_{j}}=\emptyset$ and $\varphi^{p}(U_{i})\cap U_{j}=\emptyset$ for $i\neq j$, where  $\overline{U_{i}}$ is the closure of $U_{i}$. The limit point property of $\left\lbrace\varphi^{pn}(x) |n\geq1\right\rbrace$ shows that there exists $N\geq 1$ such that $\varphi^{pn}(x) \in U_{i=1}^{k}U_{i}$ for $n\geq N$. If $\varphi^{pN}(x)\in U_{i_0}$ for some $i_0$, then, by the choice of $U_i's$, $\varphi^{pn}(x)\in U_{i_0}$ for all $n\geq N$. Hence
	the  sequence $\left\lbrace\varphi^{pn}(x) |n\geq1\right\rbrace$ has a unique limit point $z_{i_0}$, thus it converges.
	
	By a similar argument, one can obtain that $\left\lbrace\varphi^{-pn}(x) |n\geq1\right\rbrace$ converges. Thus $x\in X_p$ and the claim is established. So for each nonzero integer $p$, we have
	$\{x\in X:\, (x,\varphi^{p}(x))\in \mathcal{R}_{\varphi}\}=\{x\in X:\, (x,\varphi^{|p|}(x))\in \mathcal{R}_{\varphi}\}\subseteq X_{|p|}$, thus $\{x\in X:\, (x,\varphi^{p}(x))\in \mathcal{R}_{\varphi}\}$ is countable. Since $X$ is infinite and  has no isolated points,  it follows that $\{x\in X:\, (x,\varphi^{p}(x))\notin \mathcal{R}_{\varphi}\}$ is dense in $X$ for each nonzero integer $p$. Consequently,  $\mathbb{Z}\curvearrowright_{\alpha} (X,\mathcal{R}_{\varphi})$ is essentially free.
	
\end{proof}

Recall that the action $G\curvearrowright_{\alpha} X $ is \emph{(topologically) transitive} if for all nonempty open set $U,V \subseteq X$, there exists an $s\in G$ such that $sU\cap V\neq \emptyset$. In this case, choose a countable basis $\{U_n:\, n=1,2,\cdots\}$ for the topology on $X$. The transitivity of $\alpha$ implies that each open subset $W_n=\cup_{g\in G}\, gU_n$ is dense in $X$.
It follows from the Baire category theorem that $\cap_{n=1}^{\infty} W_n$ is dense in $X$. Thus the set of points in $X$ with dense orbit is dense.

\begin{proposition}\quad
	Let $\mathcal{R}_{\alpha}$ be the local cnjugacy relation from an expansive and transitive action $G\curvearrowright_{\alpha} X$. Assume that $X$ is infinite and has no isolated points and $G$ is an abelian group such that every subgroup generated by $g$ ($g\neq e$) has finite index in $G$. Then $G \curvearrowright_{\alpha} (X,\mathcal{R}_{\alpha})$ is essentially free.
\end{proposition}

\begin{proof}\quad
	Given $g\in G$, $g\neq e$, let $H^g$ be the subgroup generated by $g$ in $G$. Then $H^g$ has finite index,  thus $H^g\curvearrowright_{\alpha|_{H^g}} X$ is expansive, where  $\alpha|_{H^g}$ is the restriction of $\alpha$ to $H^g$. So the set $F_{H^g}(\alpha):\,=\{x\in X:\, hx=x, h\in H^g\}$ is finite. From hypothesis, $X$ is uncountable.

	If we take an enumeration $s_1, s_2, \cdots$ of the elements of $G$, then $gs_1, gs_2, \cdots $ is also an enumeration of the elements of $G$. Let $x\in X$ with $(x,gx)\in \mathcal{R}_{\alpha}$. Assume that $z$ is a limit point of   the sequence $\left\lbrace gs_nx\; :  n=1,2,\cdots \right\rbrace $. By a similar argument to the above theorem, one can see that $gz=z$, thus $z\in F_{H^g}(\alpha)$.

	Assume that there exists a $g\in G$, $g\neq e$ such that the interior of $\left\lbrace x\in X,\;(x,gx)\in \mathcal{R}_{\alpha} \right\rbrace $ is non-empty. Then the transitivity of $G\curvearrowright_{\alpha} X$ implies that there exists a point $x\in X$ with $(x,gx)\in \mathcal{R}_{\alpha}$ and having dense orbit, i.e., $\{gs_nx:\, n=1, 2, \cdots\}$ is dense in $X$. From the second paragraph, every limit point of $\left\lbrace gs_nx\; :  n=1,2,\cdots \right\rbrace $ is contained in the finite set $F_{H^g}(\alpha)$. Thus the closure of $\left\lbrace gs_nx\; :  n=1,2,\cdots \right\rbrace $ in $X$ is countable, which contradicts the fact that X is  uncountable. Consequently, $G \curvearrowright_{\alpha} (X,\mathcal{R}_{\alpha})$ is essentially free.
\end{proof}

\section{Expansive automorphism actions on compact groups}

Let $X$ be a compact metrizable group with an invariant compatible metric $d$, i.e., $d(xy,xz)=d(yx,zx)=d(y,z)$ for $x,y,z\in X$. Assume that $G\curvearrowright_{\alpha} X$ is an expansive automorphism system in the sense that it is expansive and each $\alpha_g$ is a continuous automorphism on $X$. Let $$\Delta_{\alpha}=\{x\in X:\, \lim\limits_{g\rightarrow \infty}d(\alpha_g(x),\alpha_g(e))=0\}$$
be the associated  {\sl homoclinic group}, which is an $\alpha$-invariant countable subgroup of $X$ in the sense that $\alpha_g(a)\in \Delta_{\alpha}$ for every $a\in \Delta_{\alpha}$ and $g\in G$ (\cite{Thomsen2010}).
Denote by $\sigma$ the left-multiplication action of $\Delta_{\alpha}$  on $X$:
$$\sigma_u(x)=ux, \,\, \mbox{ for $u\in \Delta_{\alpha}$ and $x\in X$},$$
and by $X\rtimes_{\sigma} \Delta_{\alpha}$ the associated transformation groupoid. Let $G\curvearrowright_{\alpha} (X,\mathcal{R})$ be the automorphism system associated to the local conjugacy equivalence relation  as in Section 5. The following facts are referred to \cite[Lemma 3.7]{Thomsen2010}.

\begin{lemma}\quad
	Let $G\curvearrowright_{\alpha} X$ be an expansive automorphism system. Then
	
	(1) two elements $x$ and $y$ in $X$ are locally conjugate, if and only if they are homoclinic, i.e., $\lim\limits_{g\rightarrow\infty}d(gx,gy)=0$, if and only if $xy^{-1}\in \Delta_{\alpha}$, if and only if $x^{-1}y\in \Delta_{\alpha}$,  if and only if $x^{-1}$ and $y^{-1}$ are locally conjugate.
	
	(2)The map $\Lambda:\, (x,y)\in \mathcal{R}\rightarrow (x,xy^{-1})\in X\rtimes_{\sigma} \Delta_{\alpha}$ is an \'{e}tale groupoid isomorphism.
	
\end{lemma}

\begin{proof}\quad		We only give a proof for (2). One can see that $\Lambda$ is an algebraic isomorphism from  $\mathcal{R}$ onto $X\rtimes_{\sigma} \Delta_{\alpha}$ with inverse map $\Lambda^{-1}$, defined by $\Lambda^{-1}(x,u)=(x,u^{-1}x)$ for $(x,u)\in X\rtimes_{\sigma} \Delta_{\alpha}$. Given $(x,y)\in \mathcal{R}$, for $S\subseteq \Delta_{\alpha}$ and an open subset $U\subseteq X$ with $x\in U$ and $xy^{-1}\in S$, we define $\gamma(z)=yx^{-1}z$ for $z\in U$.  Then $(U,\gamma(U),\gamma)$ is a local conjugacy from $x$ to $y$, and $\Lambda(\{(z, \gamma(z)):\,\, z\in U\})\subseteq U\times S$, thus $\Lambda$ is continuous at $(x,y)$. By a similar way, we show that $\Lambda^{-1}$ is continuous, thus $\Lambda$ is a homeomorphism.
	
\end{proof}

\begin{definition}\quad		
	Let $\Gamma=\Delta_{\alpha}\rtimes G$ be the semi-direct product of $\Delta_{\alpha}$ by $G$. Define the action $\widetilde{\alpha}$ of $\Gamma$ on $X$ as follows. For $(a,g)\in \Gamma$ and $x\in X$,
	$$\widetilde{\alpha}_{(a,g)}(x)=a\alpha_{g}(x).$$
\end{definition}

One can check that $\Gamma\curvearrowright_{\widetilde{\alpha}} X$ is an expansive affine system. Remark that $\Delta_{\alpha}$ and $G$ can be contained in $\Gamma$ as subgroups by identifying $a\in \Delta_{\alpha}$ with $(a,e)\in \Gamma$, and $g\in G$ with $(e,g)\in \Gamma$, thus the restrictions of $\widetilde{\alpha}$ to $\Delta_{\alpha}$ and $G$ are the same as $\sigma$ and $\alpha$, respectively.
Hence the transformation groupoid $X\rtimes_{\widetilde{\alpha}} \Gamma$  contains $X\rtimes_{\sigma} \Delta_{\alpha}$ and $X\rtimes_{\alpha} G$ as open subgroupoids.

\begin{proposition}\quad
	The map $\Lambda:(x,g,y)\longmapsto (x,(x\alpha_{g}(y^{-1}),g)) $ is an isomorphism of $\mathcal{R}\rtimes_{\alpha} G$ onto $X\rtimes_{\widetilde{\alpha}} \Gamma$ as \'{e}tale groupoids.  Moreover, $\Lambda(\mathcal{R})=X\rtimes_{\sigma}\Delta_{\alpha}$, and $\Lambda(X\rtimes_{\alpha} G)=X\rtimes_{\widetilde{\alpha}} G$.
\end{proposition}	

\begin{proof}\quad
	From Lemma 6.1, $\Lambda$ is well-defined and injective. For each $(x,(a,g))$ in $X\rtimes_{\widetilde{\alpha}} \Gamma$, we have $(x,g,\alpha_{g^{-1}}(a^{-1}x))\in \mathcal{R}\rtimes_{\alpha} G$ and
	$\Lambda(x,g,\alpha_{g^{-1}}(a^{-1}x))=(x,(a,g))$, thus $\Lambda$ is bijective. For $(x,g, y),(u,h,z)$ in $\mathcal{R}\rtimes_{\alpha} G$, we have $(x,g, y)$ and $(u,h,z)$ are composable in $\mathcal{R}\rtimes_{\alpha} G$, if and only if $u=y$, if and only if $\Lambda(x,g,y)$ and $\Lambda(y,h,z)$ are composable in $X\rtimes_{\widetilde{\alpha}} \Gamma$. Moreover,
	$$\begin{array}{lll}
		\Lambda(x,g,y)\Lambda(y,h,z) &=& (x,(x\alpha_{g}(y)^{-1},g))(y,(y\alpha_{h}(z)^{-1},h))\\
		&=& (x,(x\alpha_{gh}(z)^{-1},gh))\\
		&=& \Lambda((x,g,y)(y,h,z)).
	\end{array}$$
	
	The continuity of $\Lambda$ can be implied by Lemma 6.1  and the canonical homeomorphism $\gamma_0$ from $\mathcal{R}\rtimes_{\alpha} G$ onto $\mathcal{R}\times G$. Hence $\Lambda$ is an \'{e}tale groupoid isomorphism.
\end{proof}

\begin{proposition} \quad\begin{enumerate}
		\item[(i)] The system $G\curvearrowright_{\alpha} X$ is topologically free, if and only if $\Gamma\curvearrowright_{\widetilde{\alpha}} X$ is topologically free, if and only if
		$\mathcal{R}\rtimes_{\alpha} G$ is topologically principal, if and only if $G\curvearrowright_{\alpha} (X, \mathcal{R})$ is essentially free.
		
		\item[(ii)] If $G$ is torsion-free and $\Delta_{\alpha}$ is dense in $X$, then $G\curvearrowright_{\alpha} X$ is topologically free.
	\end{enumerate}
	
\end{proposition}

\begin{proof} \quad (i)\quad It follows from \cite[Corollary 2.3]{Li2018}, Lemma 4.6 and Proposition 6.3 that we only need to show that the topological freeness for $\alpha$ and $\widetilde{\alpha}$ is consistent. Since $G$ can be embedded into  $\Gamma$ as a subgroup and the restriction of $\widetilde{\alpha}$ to $G$ is the same as the action $\alpha$, the topological freeness of $\widetilde{\alpha}$ implies that of $\alpha$.
	
	To see the contrary, it is sufficient to show that, for arbitrary $(e,e)\neq (a,g)\in\Gamma$ and non-empty open subset $U$ of $X$, there exists $x\in U$ such that $a\alpha_g(x)\neq x$.
	
	In fact, since the restriction of $\widetilde{\alpha}$ to $\Delta_{\alpha}$ is free, we can assume that $g\neq e$ and $a\neq e$. Clearly, we can also assume that there exists $y\in U$ such that $a\alpha_g(y)=y$. The topologically freeness of $\alpha$ implies there is $z\in y^{-1}U$ such that $\alpha_g(z)\neq z$. Let $z=y^{-1}x$ for $x\in U$. Then $a\alpha_g(x)\neq x$.

	(ii)\quad Given $g\in G$, assume there exists an open subset $U$ of $X$ such that $\alpha_g(z)=z$ for every $z\in U$. We can let $e\notin U$. Since $\Delta_{\alpha}$ is dense in $X$, there is $x_0\in U\cap \Delta_{\alpha}$, thus
	$\lim\limits_{h\rightarrow \infty}d(\alpha_h(x_0),e)=0$. If $g\neq e$, then, from the torsion-freeness of $G$, the set $\{g^n:\, n\in \mathbb{Z}\}$ is infinite, we have $\lim\limits_{n\rightarrow \infty}d(\alpha_{g^n}(x_0),e)=0$, which contradicts the fact $x_0\neq e$ and $\alpha_{g^n}(x_0)=x_0$ for all $n\in \mathbb{Z}$. Consequently, $g=e$, thus $\alpha$ is topologically free.

\end{proof}

Recall that two automorphism systems $G\curvearrowright _{\alpha} X$ and $H\curvearrowright_{\beta} Y$ on compact metrizable groups are said to be \emph{algebraically conjugate} if there exist  a continuous isomorphism $\varphi:\, X\rightarrow Y$ and an isomorphism $\rho:\,G\rightarrow H$ such that $\varphi(\alpha_g(x))=\beta_{\rho(g)}(\varphi(x))$ for $g\in G$ and $x\in X$. Form \cite{Bhattacharya2000}, when $X$ and $Y$ are abelian, two notions of algebraical conjugacy and conjugacy for automorphism systems are consistent. In the following we have a similar result for automorphism actions on nonabelian groups.

\begin{proposition}\quad Let $G\curvearrowright_{\alpha} (X,\mathcal{R})$ and $H\curvearrowright_{\beta} (Y, \mathcal{S})$ be two automorphism systems on local conjugacy relations from topologically free, expansive automorphism actions on compact and connected metrizable groups $X$ and $Y$, respectively. Then the following statements are equivalent:
	\begin{enumerate}
		\item[(i)] $G\curvearrowright_{\alpha} (X,\mathcal{R})$ and $H\curvearrowright_{\beta} (Y, \mathcal{S})$ are conjugate;
		\item[(ii)] $G\curvearrowright_{\alpha} (X,\mathcal{R})\,\, \sim_{wcoe}\,\, H\curvearrowright_{\beta} (Y, \mathcal{S})$;
		\item[(iii)] $G\curvearrowright_{\alpha} X$ and $H\curvearrowright_{\beta} Y$ are continuously orbit equivalent;
		\item[(iv)] $G\curvearrowright_{\alpha} X$ and $H\curvearrowright_{\beta} Y$ are conjugate.
	\end{enumerate}
	
	Moreover, if $\Delta_{\alpha}$ is dense in $X$, then the above conditions are equivalent to the following statement.
	\begin{enumerate}
		\item[(v)] $G\curvearrowright_{\alpha} X$ and $H\curvearrowright_{\beta} Y$ are algebraically conjugate.
	\end{enumerate}

\end{proposition}

\begin{proof}\quad Since $X$ and $Y$ are connected, the continuous orbit equivalence and conjugacy of $G\curvearrowright_{\alpha} X$ and $H\curvearrowright_{\beta} Y$ are consistent. To complete the proof, we only need to  prove that $(ii)\Rightarrow (iv)$ and $(ii)\Rightarrow (v)$ when $\Delta_{\alpha}$ is dense in $X$. From Corollary 4.16 and Proposition 6.3, there is an \'{e}tale groupoid isomorphism $\Lambda:\, X\rtimes_{\widetilde{\alpha}} (\Delta_{\alpha}\rtimes G)\rightarrow Y\rtimes_{\widetilde{\beta}} (\Delta_{\beta}\rtimes H)$ such that $\Lambda(X\rtimes_{\sigma} \Delta_{\alpha})= Y\rtimes_{\sigma'} \Delta_{\beta}$, where $\sigma$ and $\sigma'$ are the left-multiplication actions, and $\widetilde{\alpha}$ and $\widetilde{\beta}$ are as in Definition 6.2.
	Since $X$ and $Y$ are connected, there are a homeomorphism $\varphi:\, X\rightarrow Y$ and a group isomorphism $\theta:\, \Delta_{\alpha}\rtimes G\rightarrow \Delta_{\beta}\rtimes H$ such that
	$$\varphi (a\alpha_g(x))=\widetilde{\beta}_{\theta(a,g)}(\varphi(x))\, \hbox{ for every $(a,g)\in \Delta_{\alpha}\rtimes G$ and $x\in X$}, \eqno{(6.1)}$$
	and $\theta(\Delta_{\alpha})=\Delta_{\beta}$, where $\Delta_{\alpha}$ and $\Delta_{\beta}$ are subgroups of the semi-direct groups as before. Define two maps $\xi:\, G\rightarrow \Delta_{\beta}$ and $\rho:\, G\rightarrow H$ by $\theta(e,g)=(\xi(g), \rho(g))$ for $g\in G$. One can check that $\rho$ is a group isomorphism by considering the inverse isomorphism $\theta^{-1}$.
	
	Letting $a=e$, the identity of $X$, in (6.1), we have $\varphi (\alpha_g(x))=\xi(g)\beta_{\rho(g)}(\varphi(x))$ for every $g\in G$ and $x\in X$. In particular, $\varphi(e)=\xi(g)\beta_{\rho(g)}(\varphi(e))$. Thus
	$\varphi (\alpha_g(x))=\varphi(e)\beta_{\rho(g)}(\varphi(e)^{-1}\varphi(x))$
	for $x\in X$ and $g\in G$. Define $\widetilde{\varphi}(x)=\varphi(e)^{-1}\varphi(x)$ for $x\in X$. Then $\widetilde{\varphi}:\, X\rightarrow Y$ is a homeomorphism and $$
	\widetilde{\varphi}(\alpha_{g}(x))=\beta_{\rho(g)}(\widetilde{\varphi}(x)) \, \hbox{ for $x\in X$ and $g\in G$}.$$
	Consequently, $G\curvearrowright_{\alpha} X$ and $H\curvearrowright_{\beta} Y$ are conjugate.
	
	Assume that $\Delta_{\alpha}$ is dense in $X$. Remark that $\theta(a,e)\in \Delta_{\beta}$, thus $\widetilde{\beta}_{\theta(a,e)}(y)=\theta(a,e)y$ for $a\in \Delta_{\alpha}$ and $y\in Y$. Letting $g=e$, the identity of $G$, and letting $x=e$, the identity of $X$,  in (6.1), one can see that $\varphi(a)=\theta(a,e)\varphi(e)$ for $a\in \Delta_{\alpha}$. Thus, by putting $g=e$ in (6.1), we have
	$\varphi(ax)=\theta(a,e)\varphi(x)=(\varphi(a)\varphi(e)^{-1}) \varphi(x)$, which implies that $\widetilde{\varphi}(ax)=\widetilde{\varphi}(a)\widetilde{\varphi}(x)$ for every $a\in \Delta_{\alpha}$ and $x\in X$. From the density of $\Delta_{\alpha}$ in $X$,  the map $\widetilde{\varphi}:\, X\rightarrow Y$ is a continuous isomorphism. So $G\curvearrowright_{\alpha} X$ and $H\curvearrowright_{\beta} Y$ are algebraically conjugate.

\end{proof}

\begin{proposition}\quad Let $G\curvearrowright_{\alpha} (X,\mathcal{R})$ be an automorphism system on local conjugacy relation from a topologically free, expansive automorphism action. Then the following statements are equivalent.
	\begin{enumerate}
		\item[(i)] $C_r^{*}(\mathcal{R})$ is simple;
		\item[(ii)] $C_r^{*}(\mathcal{R})$ has a unique tracial state;
		\item[(iii)] $\Delta_{\alpha}$ is dense;
		\item[(iv)] $C_r^*(\mathcal{R}\rtimes_{\alpha} G)$ is simple;
		\item[(v)] $C_r^*(\mathcal{R}\rtimes_{\alpha} G)$ has a unique tracial state.
	\end{enumerate}
\end{proposition}

\begin{proof}\quad For the equivalence of (i), (ii) and (iii), we refer to  \cite[Corollary 3.9]{Thomsen2010}. From Proposition 6.3, $C_r^*(\mathcal{R}\rtimes_{\alpha} G)$ is isomorphic to $C(X)\rtimes_{r,\widetilde{\alpha}} \Gamma$, thus they have the same simplicity and the uniqueness of tracial states.  From Proposition 6.4 and \cite{Li2018}, $X\rtimes_{\widetilde{\alpha}} \Gamma$ is topologically principal, thus there is  a one-to-one correspondence between the family of ideals of $C(X)\rtimes_{r,\widetilde{\alpha}} \Gamma$ and that of $\widetilde{\alpha}$-invariant open subsets of $X$ (\cite{Renault1980}).
	
	Assume (iii) holds. Since each non-empty $\widetilde{\alpha}$-invariant open subset $U$ in $X$ is invariant by the left-multiplicative by elements in $\Delta_{\alpha}$, we have $U=X$. Hence $C(X)\rtimes_{r,\widetilde{\alpha}} \Gamma$ is simple, thus (iv) holds.  On the contrary, if (iv) holds, then $C(X)\rtimes_{r,\widetilde{\alpha}} \Gamma$ is simple, which leads to the fact that the $\widetilde{\alpha}$-invariant open $X\setminus \overline{\Delta_{\alpha}}$ of $X$ is empty, where $\overline{\Delta_{\alpha}}$ is the closure of $\Delta_{\alpha}$ in $X$. Thus $\overline{\Delta_{\alpha}}=X$, i.e., (iii) holds.

	For the implication $(v) \Rightarrow (iii)$, assume that $C_r^*(\mathcal{R}\rtimes_{\alpha} G)$, thus $C(X)\rtimes_{r,\widetilde{\alpha}} \Gamma$, has a unique tracial state. If $\Delta_{\alpha}$ is not dense in $X$, then the Haar measure $\nu$ on $\overline{\Delta_{\alpha}}$ extends a Borel probability measure $\widehat{\nu}$ on $X$ different from the Haar measure $\mu_0$ on $X$. Since $\mu_0$ is invariant under the actions $\sigma$, $\alpha$ and $\widetilde{\alpha}$, for a Borel subset $E$ of $X$ and $(a,g)\in\Gamma$, we have $\widehat{\nu}(\widetilde{\alpha}_{(a,g)}(E))=\nu(\widetilde{\alpha}_{(a,g)}(E\cap \overline{\Delta_{\alpha}}))=\widehat{\nu}(E)$, thus $\widehat{\nu}$ is $\widetilde{\alpha}$-invariant. The probability measures $\widehat{\nu}$ and $\mu_0$  produce two different tracial states on $C(X)\rtimes_{r,\widetilde{\alpha}} \Gamma$, which is a contradiction.

	For the implication $(iii) \Rightarrow (v)$, assume that $\Delta_{\alpha}$ is dense in $X$. Then the Haar measure $\mu_0$ on $X$ is the unique $\widetilde{\alpha}$-invariant Borel probability measure on $X$. From \cite[Proposition 3.2.4]{Tomiyama1987}, $C(X)\rtimes_{r,\widetilde{\alpha}} \Gamma$, and thus $C_r^*(\mathcal{R}\rtimes_{\alpha} G)$ has a unique tracial state.
\end{proof}

\begin{example}[Hyperbolic toral automorphisms ]\quad
	For $n\geq 2$, we consider an expansive  $\mathbb{Z}$-action on  the $n$-dimensional torus $\mathbb{R}^n/\mathbb{Z}^n$ generated by a single hyperbolic toral automorphism $\alpha$. Let  $\pi:\,\, \mathbb{R}^n\rightarrow \mathbb{R}^n/\mathbb{Z}^n$ be the usual quotient map. Recall that  $\mathbb{R}^n/\mathbb{Z}^n$ is a compact and connected  additive group under the following metric compatible with the quotient topology:
	$$d(\pi(x),\pi(y))=\inf\limits_{z\in \mathbb{Z}^n}\|x-y-z\|, \hbox{ for $x,y\in \mathbb{R}^n$},$$
	where $\|\cdot\|$ is the Euclidean norm on $\mathbb{R}^n$. The elements in $\mathbb{R}^n$ are denoted by column vectors or row vectors.
	
	Let $A$ be the hyperbolic matrix in $GL(n,\mathbb{Z})$ with $det(A)=\pm 1$ and having no eigenvalues of modules $1$, such that
	$$\alpha(\pi(x))=\pi(Ax)\hbox{ for $x\in \mathbb{R}^n$}.$$
	Then $\mathbb{R}^n=E^s\oplus E^u$,  where $E^s=\{x\in \mathbb{R}^n:\, \lim\limits_{k\rightarrow +\infty}A^kx=0\}$ and $E^u=\{w\in \mathbb{R}^n:\, \lim\limits_{k\rightarrow +\infty}A^{-k}w=0\}$ are two invariant subspaces of the linear map on $\mathbb{R}^n$ determined by  $A$. Remark that $E^s\cap \mathbb{Z}^n=\{0\}$, $E^u\cap \mathbb{Z}^n=\{0\}$, and both subgroups $\pi(E^s)$ and $\pi(E^u)$, as well as
	the homoclinic group $\Delta_{\alpha}=\pi(E^s)\cap \pi(E^u)$ induced by $\alpha$,
	are dense in $\mathbb{R}^n/\mathbb{Z}^n$. Moreover, the system $\mathbb{Z}\curvearrowright_{\alpha} \mathbb{R}^n/\mathbb{Z}^n$ generated by $\alpha$  is topologically free (\cite{Lind1999}).
	
	Each ${\bf m}\in \mathbb{Z}^n$ has the unique decomposition ${\bf m}={\bf m}_s-{\bf m}_u\in E^s\oplus E^u$. Then the map $\theta: \mathbb{Z}^n\rightarrow \Delta_{\alpha}$ by $\theta(\bf{m})=\pi({\bf m}_s)$ ($=\pi({\bf m}_u)$) is a group isomorphism. As before, we let $\sigma$ be the translation action of $\Delta_{\alpha}$ on $\mathbb{R}^n/\mathbb{Z}^n$: $$\sigma_u(x)=u+x\,\, \hbox{ for $u\in \Delta_{\alpha}, x\in \mathbb{R}^n/\mathbb{Z}^n$.}$$
	Let $\tau$ be the action of $\mathbb{Z}^n$ on $\mathbb{R}^n/\mathbb{Z}^n$ by homeomorphisms: $$\tau_{\bf n}(x)=\theta({\bf n})+x\,\, \hbox{ for ${\bf n}\in \mathbb{Z}^n, x\in \mathbb{R}^n/\mathbb{Z}^n$.}$$
	Then $ \mathbb{Z}^n \curvearrowright_{\tau}\mathbb{R}^n/\mathbb{Z}^n$ and $\Delta_{\alpha} \curvearrowright_{\sigma}\mathbb{R}^n/\mathbb{Z}^n$ are conjugate.
	
	Denote by $\mathbb{Z}^n\rtimes \mathbb{Z}$ the semi-direct product of $\mathbb{Z}^n$ by the automorphism given by $A$: ${\bf m}\in \mathbb{Z}^n\rightarrow A{\bf m}\in \mathbb{Z}^n$. Let $\gamma$ be the action of $\mathbb{Z}^n\rtimes \mathbb{Z}$ on $\mathbb{R}^n/\mathbb{Z}^n$:
	$$\gamma_{({\bf m},k)}(x)=\theta({\bf m})+\alpha^k(x)\,\, \hbox{ for $({\bf m},k)\in \mathbb{Z}^n\rtimes \mathbb{Z}$ and $x\in \mathbb{R}^n/\mathbb{Z}^n$}.$$
	So $\mathbb{Z}^n\rtimes \mathbb{Z} \curvearrowright_{\gamma} \mathbb{R}^n/\mathbb{Z}^n$ and $\Delta_{\alpha}\rtimes \mathbb{Z}\curvearrowright_{\widetilde{ \alpha}}  \mathbb{R}^n/\mathbb{Z}^n $ are conjugate, where $\widetilde{\alpha}$ is given by Definition 6.2.

	We consider the multiplicative coordinate system on the $n$-dimensional torus by $\mathbb{T}^n:\, =\{(z_1,z_2,\cdots,z_n):\, z_i\in \mathbb{C}, |z_i|=1, \hbox{ for $1\leq i\leq n$}\}$. The correspondence $$\varphi:\, [(x_1,x_2,\cdots,x_n)]\in \mathbb{R}^n/\mathbb{Z}^n\rightarrow (e^{2\pi \mathrm{i}x_{1}},e^{2\pi \mathrm{i}x_{2}},\cdots,e^{2\pi \mathrm{i}x_{n}})\in \mathbb{T}^n$$ is an isomorphism between two representations, where
	$[(x_1,x_2,\cdots,x_n)]\,=\pi(x_1,x_2,\cdots,x_n)$ for $(x_1,x_2,\cdots,x_n)\in \mathbb{R}^n$. Using this coordinate system, we can rewrite the above toral automorphism $\alpha$ and the actions $\tau$, $\gamma$ as follows. Let $A=(a_{ij})$ and $A^{-1}=(b_{ij})$. Define the automorphism $\beta$ of $\mathbb{T}^n$ by
	$$\beta(z_1,z_2,\cdots,z_n)=(z_1^{a_{11}}z_2^{a_{12}}\cdots z_n^{a_{1n}}, z_1^{a_{21}}z_2^{a_{22}}\cdots z_n^{a_{2n}},\cdots, z_1^{a_{n1}}z_2^{a_{n2}}\cdots z_n^{a_{nn}})$$
	for $(z_1,z_2,\cdots,z_n)\in \mathbb{T}^n$,
	the rotation action $\rho$ of $\mathbb{Z}^n$ on $\mathbb{T}^n$ by
	$$ \rho_{\bf m} (v)=\varphi(\theta({\bf m}))v\, \hbox{ for ${\bf m}\in \mathbb{Z}^n$ and $v\in \mathbb{T}^n$},$$
	and the action $\widetilde{\gamma}$ of $\mathbb{Z}^n\rtimes \mathbb{Z}$ on $\mathbb{T}^n$ by
	$$\widetilde{\gamma}_{({\bf m},k)}(v)=\varphi(\theta({\bf m}))\beta^k(v)\,\, \hbox{ for $({\bf m},k)\in \mathbb{Z}^n\rtimes \mathbb{Z}$ and $v\in \mathbb{T}^n$}.$$
	Then $\mathbb{Z}^n \curvearrowright _{\tau}\mathbb{R}^n/\mathbb{Z}^n$ and $\mathbb{Z}^n\curvearrowright_{\rho}\mathbb{T}^n $ are conjugate, and $ \mathbb{Z}^n\rtimes \mathbb{Z}\curvearrowright_{\widetilde{\gamma}}\mathbb{T}^n $ and $ \mathbb{Z}^n\rtimes \mathbb{Z}\curvearrowright_{\gamma}\mathbb{R}^n/\mathbb{Z}^n $ are  conjugate.
	
	From Lemma 6.1, Proposition 6.3 and the above, the local conjugacy relation $\mathcal{R}$ given by $\alpha$  and the associated semi-direct product $\mathcal{R}\rtimes_{\alpha} \mathbb{Z}$ are isomorphic to  the transformation groupoids $\mathbb{T}^n\rtimes_{\rho} \mathbb{Z}^n$ and $\mathbb{T}^n\rtimes_{\widetilde{\gamma}} (\mathbb{Z}^n\rtimes \mathbb{Z})$, respectively.
	
	We still denote by $\rho$ the automorphism action of $\mathbb{Z}^n$ on $C(\mathbb{T}^n)$ induced by the system $\mathbb{Z}^n\curvearrowright_{\rho}\mathbb{T}^n$: $$\rho_{\bf m}(f)(v)=f(\varphi(\theta({\bf m}))^{-1}v)$$ for ${\bf m}\in \mathbb{Z}^n$, $f\in C(\mathbb{T}^n)$ and $v\in \mathbb{T}^n$. Let ${\bf e}_k, 1\leq k\leq n,$ be the canonical basis of $\mathbb{Z}^n$ and $\theta(\mathit{e_{k}})=[(\theta_{k1},\cdots,\theta_{kn})]\in \Delta_{\alpha}$, where $\theta_{kj}\in[0,1]$.
	Let $U_j$, $1\leq j\leq n$, be the unitaries in $C(\mathbb{T}^n)$ defined by $U_j(z_1,\cdots,z_n)=z_j$ for $(z_1,\cdots,z_n)\in \mathbb{T}^n$, and let $V_k$, $1\leq k\leq n$, be the unitaries implementing the $C^*$-automorphism $\rho_{{\bf e}_k}$ on $C(\mathbb{T}^n)$. One can check that $$ U_jU_k=U_kU_j,\, V_jV_k=V_kV_j,\,\, U_jV_k=e^{2\pi i\theta_{kj}}V_kU_j \eqno{(6.2)}$$
	for $1\leq j,k\leq n$. From Proposition 6.4 and 6.6, $C_r^*(\mathcal{R})$, thus $C(\mathbb{T}^n)\rtimes_{\rho} \mathbb{Z}^n$ are simple and have unique tracial states. Hence $C(\mathbb{T}^n)\rtimes_{\rho} \mathbb{Z}^n$ is generated by $U_j,\, V_j$, $1\leq j\leq n$, thus is the $2n$-dimensional noncommutative torus $A_{\Theta}$ for a $2n\times 2n$ real skew-symmetric matrix $\Theta=(\widetilde{\theta}_{kl})$ defined by
	$\widetilde{\theta}_{kl}=0$ for $1\leq k,l\leq n$ or $n+1\leq k,l\leq 2n$, $\widetilde{\theta}_{kl}=\theta_{k(l-n)}$ for $1\leq k\leq n$ and $n+1\leq l\leq 2n$, and $\widetilde{\theta}_{kl}=-\theta_{l(k-n)}$ for $n+1\leq k\leq 2n$ and $1\leq l\leq n$ (\cite{Rieffel1990}). From \cite{Phillips2006}, $C_r^*(\mathcal{R})$ is an AT-algebra with real rank zero and the range of the unique tracial state acting on $K_0(C_r^*(\mathcal{R}))$ is an isomorphism invariant.

	Similarly, we also denote by $\widetilde{\gamma}$ the automorphism action of $\mathbb{Z}^n\rtimes \mathbb{Z}$ on $C(\mathbb{T}^n)$ induced by the system $ \mathbb{Z}^n\rtimes \mathbb{Z}\curvearrowright_{\widetilde{\gamma}}\mathbb{T}^n $:
	$$\widetilde{\gamma}_{({\bf m},k)}(f)(v)=f((\varphi\theta(A^{-k}{\bf m}))^{-1}\cdot \beta^{-k}(v))$$ for $({\bf m},k)\in \mathbb{Z}^n\rtimes \mathbb{Z}$, $f\in C(\mathbb{T}^n)$ and $v\in \mathbb{T}^n$. Let $U_j$, $j=1,2,\cdots,n$, be the generating set of $C(\mathbb{T}^n)$ as   above, and let $V_j'$, $j=1,2,\cdots,n$ and $W$ be the unitaries implementing the automorphisms $\widetilde{\gamma}_{({\bf e}_j,0)}$ and $\widetilde{\gamma}_{(0,1)}$ associated to the generating set $({\bf e}_j,0)$ for $1\leq j\leq n$ and $(0,1)$ of $\mathbb{Z}^n\rtimes \mathbb{Z}$. Then we have
	$$
	\begin{array}{ll}
		U_jU_k=U_kU_j,\, V_j'V_k'=V_k'V_j',\,\, U_jV_k'=e^{2\pi i\theta_{kj}}V_k'U_j\\
		\\
		WU_jW^*=\prod\limits_{l=1}^nU_l^{b_{jl}},\,\,\, WV_j'W^*=\prod\limits_{l=1}^nV_l'\,^{a_{lj}}
	\end{array} \eqno{(6.3)}$$
	for $1\leq j,k\leq n$. Since $C(\mathbb{T}^n)\rtimes_{\widetilde{\gamma}} (\mathbb{Z}^n\rtimes \mathbb{Z})$ is simple from Proposition 6.4 and 6.6, it is generated by the unitaries $U_j, V_j'$, $1\leq j\leq n$, and $W$ satisfying the above relations.

	From the above arguments and \cite[Theorem 3.33 and 3.36]{Thomsen2010}, we have the following results which generalize \cite[Theorem 2.9 and Proposition 6.1]{Matsumoto20212}.
	
	\begin{proposition}\quad
		Let $\alpha$ be a hyperbolic toral automorphism on $\mathbb{R}^n/\mathbb{Z}^n$ defined by a hyperbolic matrix $A$. Let $\mathcal{R}$ be the local conjugacy relation induced by $\alpha$.  Then
		\begin{enumerate}
			\item[(1)] $C_r^*(\mathcal{R})$ is generated by the unitaries $U_j, \, V_j$, $1\leq j\leq n$, satisfying the relations $(6.2)$, thus is isomorphic to a simple $2n$-dimensional noncommutative torus and is an AT-algebra with real rank zero.
			\item[(2)] $C_r^*(\mathcal{R}\rtimes_{\alpha} \mathbb{Z})$ is generated by unitaries $U_j,\, V_j'$, $1\leq j\leq n$, and $W$ satisfying the relations $(6.3)$.
		\end{enumerate}
		
		Moreover, two hyperbolic toral automorphisms  on $\mathbb{R}^n/\mathbb{Z}^n$ are flip conjugate if and only if the $\mathbb{Z}$-actions they generates are continuously orbit equivalent up to the associated local conjugacy relations.
	\end{proposition}

\end{example}

\subsection*{Acknowledgements} This work was supported by the NSF of China (Grant No. 12271469, 11771379, 11971419).


\end{document}